\documentclass[oneside,english]{sn-jnl}
\usepackage[T1]{fontenc}
\usepackage[utf8]{inputenc}
\usepackage[active]{srcltx}
\usepackage{textcomp}
\usepackage{amsmath}
\usepackage{amssymb}
\usepackage{graphicx}

\makeatletter

\newcommand{\lyxmathsym}[1]{\ifmmode\begingroup\def\b@ld{bold}
  \text{\ifx\math@version\b@ld\bfseries\fi#1}\endgroup\else#1\fi}

\makeatother

\usepackage{babel}
\begin{document}
\title[Preserving Superconvergence of Spectral Elements]{Preserving
Superconvergence of Spectral Elements for Curved Domains via $h$
and $p$-Geometric Refinement}

\author[1]{\fnm{Jacob} \sur{Jones}} \email{jacob.jones@stonybrook.edu}
\author[2]{\fnm{Rebecca} \sur{Conley}} \email{rconley@saintpeters.edu}
\author[1]{\fnm{Xiangmin} \sur{Jiao}}\email{xiangmin.jiao@stonybrook.edu}
\affil[1]{\orgdiv{Dept. of Applied Mathematics \& Statistics and Institute for Advanced Computational Science}, \orgname{Stony Brook University}, \orgaddress{ \city{Stony Brook}, \state{NY}, \country{USA}}}
\affil[2]{\orgdiv{Dept. of Mathematics \& Statistics}, \orgname{Saint Peter's University}, \orgaddress{\city{Jersey City}, \state{NJ}, \country{USA}}}
\abstract{Spectral element methods (SEM), which are extensions of finite element
methods (FEM), are important emerging techniques for solving partial
differential equations in physics and engineering. SEM can potentially
deliver better accuracy due to the potential superconvergence for
well-shaped tensor-product elements. However, for complex geometries,
the accuracy of SEM often degrades due to a combination of geometric
inaccuracies near curved boundaries and the loss of superconvergence
with simplicial or non-tensor-product elements. We propose to overcome
the first issue by using $h$- and $p$-geometric refinement, to
refine the mesh near high-curvature regions and increase the degree
of geometric basis functions, respectively. We show that when using
mixed-meshes with tensor-product elements in the interior of the domain,
curvature-based geometric refinement near boundaries can improve the
accuracy of the interior elements by reducing pollution errors and
preserving the superconvergence. To overcome the second issue, we
apply a post-processing technique to recover the accuracy near the
curved boundaries by using the adaptive extended stencil finite element
method (AES-FEM). The combination of curvature-based geometric refinement
and accurate post-processing delivers an effective and easier-to-implement
alternative to other methods based on exact geometries. We demonstrate
our techniques by solving the convection-diffusion equation in 2D
and show one to two orders of magnitude of improvement in the solution
accuracy, even when the elements are poorly shaped near boundaries.}

\maketitle
\maketitle

\section{Introduction}

Spectral element methods (SEM) are extensions of finite element methods
(FEM) that use Gauss-Lobatto or similar nodes instead of equidistant
nodes for high-order elements \cite{karniadakis2005spectral}. Mathematically,
SEM with well-shaped tensor-product elements can deliver higher accuracy
than equidistant FEM due to potential superconvergence \cite{kvrivzek1987superconvergence,Chen1995,Zlamal1977,Zlamal1978},
in that the $\ell_{2}$-norm error (i.e., nodal solution error) may
converge faster (at $(p+2)$nd order with degree-$p$ elements \cite{kvrivzek1987superconvergence,Chen1995,Zlamal1977,Zlamal1978})
than the standard “optimal” rate of the $L_{2}$ norm (at $(p+1)$st
order) of FEM. However, some significant challenges remain for domains
with curved boundaries, limiting the advantages of SEM for some real-world
applications. In particular, non-tensor-product elements are typically
used near boundaries, leading to the loss of superconvergence of the
overall solution. In addition, for curved domains, isoparametric elements
are commonly used, for which the geometry is approximated to the same
order as the solution space. There has been significant evidence in
isogeometric analysis \cite{hughes2005isogeometric,scott2011isogeometric}
and NURBS-enhanced FEM (NEFEM) \cite{sevilla2008nurbs,sevilla2011comparison}
that a more accurate boundary representation can significantly improve
the overall accuracy of the solutions. Similar observations have been
made in the context of discontinuous Galerkin methods \cite{botti2018assessment}.
However, to the best of our knowledge, there has not been a systematic
study on the issue in the context of continuous Galerkin formulations
with high-order finite or spectral elements. 

In this work, we propose a novel approach to improve the overall accuracy
and preserve the superconvergence of SEM over curved domains. Our
approach has two distinct components. First, we develop a mesh-generation
procedure using rectangular tensor-product elements in the interior
of the domain but geometrically refined non-tensor-product elements
near curved boundaries. The geometric refinement includes curvature-based
edge-length controls and optionally superparametric elements, which
we refer to as $h$- and $p$-geometric refinements ($h$-GR and $p$-GR),
respectively. We show that the refinements can improve the accuracy
of the geometric representations and, in turn, reduce the pollution
errors of the tensor-product elements in the interior of the domain.
Second, we propose a post-processing phase to improve the accuracy
of the near-boundary elements using the adaptive extended stencil
finite element method (AES-FEM) \cite{conley2016overcoming,conley2020hybrid},
so that the accuracy of non-tensor-product elements can match those
of the superconvergent spectral elements. Overall, we show that our
approach can deliver one to two orders of improvement in solution
accuracy. 

The contributions of the work can be summarized as follows. First,
to the best of our knowledge, it is the first work to demonstrate
the superconvergence of SEM over curved domains with both Dirichlet
and Neumann boundary conditions. This result is significant in that
theoretical analysis of superconvergence is limited to only Dirichlet
boundary conditions over rectangular domains \cite{kvrivzek1987superconvergence,Chen1995,Zlamal1977,Zlamal1978}.
Second, we show that contrary to the conventional wisdom of isoparametric
elements \cite{strangfem}, superparametric elements can play an important
role in improving the accuracy of nodal solutions, especially when
combined with $h$-GR. Third, we show that the impact of geometric
refinements may be orders of magnitude more significant than the element
shapes for FEM and SEM. These results suggest that we should focus
on geometric accuracy instead of element shapes for high-order elements.
 Fourth, our overall workflow offers a promising alternative to other
discretization methods using exact geometries, such as NEFEM \cite{sevilla2008nurbs,sevilla2011comparison}.
Our approach is easier to implement because its mesh generation and
adaptation procedure is more consistent with standard FEM and SEM.
Its post-processing step can be incorporated into other superconvergence
post-processing steps with a relatively small computational cost.

The remainder of the paper is organized as follows. In section~\ref{sec:Background},
we present some background information on SEM and related pre- and
post-processing techniques. In section~\ref{sec:sem-curved-domains},
we discuss the reasons for the loss of superconvergence of SEM in
the presence of curved boundaries and some potential remedies. In
section~\ref{sec:Mesh-Generation-Workflow}, we describe our mesh-generation
procedure, focusing on $h$- and $p$-GR and the placement of high-order
nodes in 2D, with a brief discussion on the extension to 3D. In section~\ref{sec:Results},
we present numerical results and demonstrate the impact of geometric
refinement on the preservation of superconvergence. Section~\ref{sec:Conclusions}
concludes the paper.

\section{Background\label{sec:Background}}

We first review FEM and SEM, focusing on their superconvergence properties
and pre- and post-processing techniques over curved domains. 

\subsection{\label{subsec:Finite-and-Spectral}Finite and Spectral Element Methods}

Mathematically, FEM and SEM are spatial discretization techniques
for partial differential equations (PDEs). For simplicity, let us
consider linear elliptic PDEs, such as the convention-diffusion equation,
which we will solve in section~\ref{sec:Results}. Let $\Omega\subset\mathbb{\mathbb{R}}^{d}$
be a bounded, piecewise smooth domain with boundary $\partial\Omega=\Gamma_{D}\cup\Gamma_{N}$,
where $\Gamma_{D}$ and $\Gamma_{N}$ denote the Dirichlet and Neumann
boundaries, respectively. Typically, $d=2$ or 3. Let $\mathcal{L}$
be a second-order elliptic linear differential operator. The associated
boundary value problem is to find a sufficiently smooth solution $u$
such that
\begin{align}
\mathcal{L}u=f & \quad\text{on }\Omega,\label{eq:BVP PDE}\\
u=u_{D} & \quad\text{on }\Gamma_{D},\label{eq:BVP Dirchlet}\\
\partial_{\boldsymbol{n}}u=g_{N} & \quad\text{on }\Gamma_{N},\label{eq:BVP Neumann}
\end{align}
where $u_{D}$ denotes the Dirichlet boundary conditions and $g_{N}$
denotes the Neumann boundary conditions. The domain $\Omega$ is tessellated
into a mesh, composed of a set of elements. FEM and SEM approximate
the solution to equations (\ref{eq:BVP PDE})--(\ref{eq:BVP Neumann})
using the method of weighted residuals \cite{finlayson2013method}.
Let $u_{h}$ denote the approximate solution, and the residual of
equation (\ref{eq:BVP PDE}) is then $\mathcal{L}u_{h}-f$. Let $\Psi=\left\{ \psi_{i}\mid1\leq i\leq n\right\} $
denote a set of test functions. A weighted residual method requires
the residual to be orthogonal to $\Psi$ over $\Omega$, that is,
\[
\int_{\Omega}\mathcal{L}u_{h}\psi_{i}\ d\Omega=\int_{\Omega}f\psi_{i}\ d\Omega,\quad\text{for }i=1,\dots,n.
\]
Traditionally, FEM uses equidistant nodes within each element and
the basis functions for $u_{h}$ are piecewise Lagrange polynomial
basis functions based on these nodes. In contrast, SEM uses Gauss-Lobatto
nodes over tensor-product elements \cite{patera1984spectral,maday1990optimal,karniadakis2005spectral},
based on the tensor-product Gauss-Lobatto quadrature rules. One of
the motivations for using Gauss-Lobatto points is the improved stability
of polynomial interpolations for very high-degree polynomials. However,
for moderate degrees, the major advantage of SEM is that the nodal
solutions naturally \emph{superconverge}, in the sense that the $\ell_{2}$-norm
error is $\mathcal{O}(h^{p+2})$ for degree-$p$ elements \cite{chen1981superconvergence},
whereas the optimal convergence rate in $L_{2}$ norm is only $\mathcal{O}(h^{p+1})$
\cite{strangfem}, where $h$ denotes an edge-length measure, assuming
the mesh is quasi-uniform. 

The superconvergence property of tensor-product SEM has been studied
extensively in the literature; we refer readers to \cite{brenner2004finite,zhu1998review}
and references therein.  In a nutshell, the superconvergence is related
to the fact that the Gauss-Lobatto polynomials are derivatives of
Legendre polynomials \cite{chen1981superconvergence}. However, this
connection with orthogonal polynomials does not generalize to non-tensor-product
elements, regardless of whether they are obtained by degenerating
tensor-product elements \cite{karniadakis2005spectral} or using some
Gauss-Lobatto-like point distributions (such as Fekete points \cite{taylor2000algorithm}).
Hence, non-tensor-product ``spectral'' elements generally do not
superconverge. Nevertheless, they are useful for compatibility with
tensor-product spectral elements in mixed-element meshes (see e.g.,
\cite{karniadakis2005spectral,komatitsch2001wave}), which are desirable
to handle complex geometries arising from engineering applications.
In this work, we use mixed-element meshes as in \cite{karniadakis2005spectral,komatitsch2001wave}.
Still, we aim to preserve the superconvergence of tensor-product elements
and, in addition, improve the accuracy of the non-tensor-product element.

\subsection{Resolution of Curved Boundaries in Finite Element Analysis (FEA)\label{subsec:Resolution-of-Curved}}

Curved boundaries have long been recognized as a significant source
of error in FEA because insufficient geometric accuracy can lead to
suboptimal convergence rates (see e.g., \cite{cheung2019optimally,luo2001influence}).
Traditionally, isoparametric elements were believed to be sufficient
to preserve the $\mathcal{O}(h^{p+1})$ accuracy in $L_{2}$ norm
\cite{strangfem}. Nevertheless, superparametric elements, in which
the degree of the geometric basis functions is greater than the degree
of the trial space (see e.g., \cite{eslami2014finite}), have been
shown to improve accuracy in the context of discontinuous Galerkin
methods \cite{zwanenburg2017necessity,bassi1997high}.  In this
work, we will show that superparametric elements can also significantly
improve accuracy for continuous Galerkin methods in $\ell_{2}$ norm
(i.e., nodal solutions), especially for Neumann boundary conditions,
thanks to the improved accuracy in normal directions. 

In recent years, some alternatives have been developed to improve
geometric accuracy. Most notably, isogeometric analysis (IGA) \cite{hughes2005isogeometric}
and NURBS-enhanced finite element method (NEFEM) \cite{sevilla2008nurbs,sevilla2011comparison}
use exact geometries. IGA is analogous to isoparametric FEM, using
the same basis functions for the geometric representation and the
solutions. Since geometric models (such as NURBS) typically use $C^{1}$
or $C^{2}$ basis functions, IGA has higher regularity of its solutions,
allowing it to resolve smoother functions more effectively in some
applications. In contrast, NEFEM is analogous to superparametric FEM
in that it uses the exact geometric representation while preserving
the Lagrange polynomial basis of the trial and test spaces of FEM.
Both IGA and NEFEM pose significant challenges in mesh generation
for complex domains. 

\subsection{Curvature-Based Mesh Generation and Adaptation}

Using superparametric elements or exact geometry increases the degree
of the geometric representations. In addition, one can also improve
geometric accuracy by adapting the edge length based on curvature.
Such an approach is relatively popular in computer graphics and visualization;
see e.g., the survey paper \cite{khan2020surface} and references
therein. For example, in \cite{dassi2014curvature}, Dassi et al.
introduced an optimization-based surface remeshing method, whether
the energy is defined based on curvatures. In \cite{mansouri2016segmentation},
Mansouri and Ebrahimnezhad used Voronoi tesselation and Lloyd's relaxation
algorithm \cite{lloyd1982least} to create curvature adaptive site
centers. These site centers are used to generate a curvature-based
mesh.  In the context of FEM, Moxey et al. \cite{moxey2015isoparametric}
generated boundary layer meshes based on curvatures.  In this work,
we use the radii of curvatures to guide the refinement of edge lengths.
We also use it in conjunction with superparametric elements to improve
geometric accuracy, to reduce the pollution errors of the spectral
elements in the interior of the domain.

\subsection{Superconvergent Post-Processing Techniques}

An important technique in this work is to post-process the solutions
near boundaries to recover higher accuracy near boundaries. To this
end, we use AES-FEM, a method based on weighted least-squares approximations
\cite{conley2016overcoming,conley2020hybrid}. We defer more detailed
discussions on AES-FEM to section~\ref{subsec:Recovery-of-Accuracy}.
In FEM, post-processing is often used to improve the accuracy of the
solutions or gradients. For example, Bramble and Schatz \cite{bramble1977higher}
developed a post-processing method to achieve superconvergence of
the solution using an averaging method, but the technique requires
the mesh to be locally translation invariant. In \cite{Babuska1984},
Babuška and Miller constructed generating functions to compute weighted
averages of solutions over the entire finite element domain, but the
method did not produce superconvergent results.  The post-processing
of gradients has enjoyed better success.  One such technique is the
superconvergent patch recovery (SPR) method, which constructs a local
polynomial fitting on an element patch and then evaluates the fittings
at Gaussian quadrature points to achieve superconvergent first derivatives
\cite{Zienkiewicz1992b}. Other examples include the polynomial preserving
recovery (PPR) method \cite{Naga2004} and recovery by equilibrium
in patches (REP) \cite{Boroomand1997,Boroomand1997a}, which also
construct local polynomial fittings. One similarity of these post-processing
techniques with AES-FEM is using some form of least-squares fittings.
However, unlike those techniques, AES-FEM uses weighted least-squares
to compute differential operators and then uses them to solve the
PDEs to achieve superconvergent solutions.

\section{\label{sec:sem-curved-domains}Spectral Elements over Curved Domains}

As discussed in section~\ref{subsec:Finite-and-Spectral}, spectral
elements can enjoy superconvergence over rectangular domains, but
curved domains pose several challenges. We discuss these challenges
and the potential remedies in this section.

\subsection{Loss and Preservation of Superconvergence\label{subsec:Loss-and-Preservation}}

The theoretical analysis of SEM's superconvergence property is generally
limited to rectangular domains with homogeneous (Dirichlet) boundary
conditions \cite{Chen1995,chen1999superconvergence,Chen2005}. The
extension of such analysis to curved domains with general boundary
conditions has fundamental challenges, and this work aims to address
these challenges.  First, the most systematic analyses, such as element
orthogonality analysis (EOA) and orthogonality correction technique
(OCT) \cite{chen1981superconvergence,Chen1995,Chen2005,chen2013highest},
rely on the properties of orthogonal polynomials. However, the orthogonality
is lost if tensor-product elements are used near curved boundaries.
Mathematically, this loss is because the Lagrange basis functions
in SEM (as in FEM) within each element, say $\phi_{i}(\boldsymbol{x})$,
are defined as
\[
\phi_{i}(\boldsymbol{x})=\phi_{i}(\boldsymbol{\xi}(\boldsymbol{x})),
\]
where $\boldsymbol{\xi}(\boldsymbol{x})$ is the inverse of the mapping
\[
\boldsymbol{x}(\boldsymbol{\xi})=\sum_{k=1}^{n_{e}}\boldsymbol{x}_{k}\varphi_{k}(\boldsymbol{\xi}).
\]
In the above, $\boldsymbol{\xi}$ denotes the local (or natural) coordinates
within an element $\tau$, and $\varphi_{k}(\boldsymbol{\xi})$ denotes
the geometric basis function corresponding to the $k$th node of $\tau$.
Although the $\phi_{i}(\boldsymbol{\xi})$ and $\varphi_{k}(\boldsymbol{\xi})$
are polynomials, the $\phi_{i}(\boldsymbol{x})$ are no longer polynomials
if $\boldsymbol{x}(\boldsymbol{\xi})$ is nonlinear, and as a result,
orthogonality is also lost. Hence, we cannot expect superconvergence
of the spectral elements near boundaries even if tensor-product spectral
elements are used. It is worth noting that this loss of orthogonality
also applies to the interior of the domain if the tensor-product elements
are skewed (i.e., no longer rectangular). Therefore, we use rectangular
spectral elements in the interior of the domain as much as possible.
Near boundaries, it suffices to use non-tensor-product elements for
more flexibility in mesh generation. These elements cannot superconverge,
and we address their accuracy issues in section~\ref{subsec:Recovery-of-Accuracy}.

Second, even if rectangular spectral elements are used in the interior
of the domain, the nodal solutions of those elements may still lose
superconvergence if the near-boundary elements are too inaccurate
because the stiffness matrix couples the nodal solutions of all the
elements. We refer to these potential errors in the spectral elements
due to inaccurate near-boundary elements as \emph{pollution errors}.
These pollution errors are especially problematic when Neumann boundary
conditions are applied to curved boundaries with isoparametric elements
since the normal directions are, in general, a factor of $1/h$ less
accurate than the nodal positions. For elliptic PDEs, these pollution
errors decay with respect to the distance to the boundary. Nevertheless,
it is important to reduce the errors of near-boundary elements as
much as possible. This analysis motivates the $h$- and $p$-GR that
we will describe in section~\ref{sec:Mesh-Generation-Workflow}.
We will verify this analysis numerically in section~\ref{sec:Results}.

\subsection{\label{subsec:Recovery-of-Accuracy}Recovery of Accuracy Near Boundaries}

As noted above, the superconvergence is lost for non-tensor-product
elements near boundaries. In other words, the nodal solutions near
boundaries will be limited to only $\mathcal{O}(h^{p+1})$. Note that
this limitation cannot be overcome using $h$-adaptivity within the
FEM framework since the non-tensor-product elements directly adjacent
to tensor-product elements cannot decrease the edge length. Neither
can it be overcome using $p$-adaptivity with non-conformal meshes
since the continuity constraints of non-conformal meshes would cause
the tensor-product elements to lose superconvergence.

To overcome this fundamental challenge, we utilize the high-order
adaptive extended stencil finite element method (AES-FEM) as described
in \cite{conley2020hybrid}. AES-FEM is a generalized weighted residual
method, of which the trial functions are generalized Lagrange polynomial
(GLP) basis functions constructed using local least squares fittings,
but the test functions are $C^{0}$ continuous Lagrange basis functions
as in FEM. As shown in \cite{conley2020hybrid}, on a sufficiently
fine quasiuniform mesh with consistent and stable GLP basis functions,
the solution of a well-posed even-degree-$p$ AES-FEM for a coercive
elliptic PDE with Dirichlet boundary conditions converges \textit{\emph{at}}\emph{
$\mathcal{O}(h^{p})$ }\textit{\emph{in}}\emph{ $\ell^{\infty}$}\textit{\emph{
norm}}\emph{, }if $p=2$ or if $p>2$ and the local support is nearly
symmetric. We refer readers to \cite{conley2020hybrid} for more details.
In this work, we apply AES-FEM as a post-processing step for near-boundary
elements by splitting the near-boundary elements into linear elements,
as we describe in section~\ref{subsec:near-boundary-aes-fem}.

\section{Curvature-Based $h$- and $hp$-Geometric Refinement\label{sec:Mesh-Generation-Workflow}}

We now describe the generation of mixed-element meshes, which is the
core of this work. We describe the mesh-generation process in two
dimensions and then outline how the process can be extended to three
dimensions.

\subsection{Initial Mesh Setup}

\begin{figure*}
\begin{minipage}[t]{0.33\textwidth}%
\begin{center}
\includegraphics[width=0.9\textwidth]{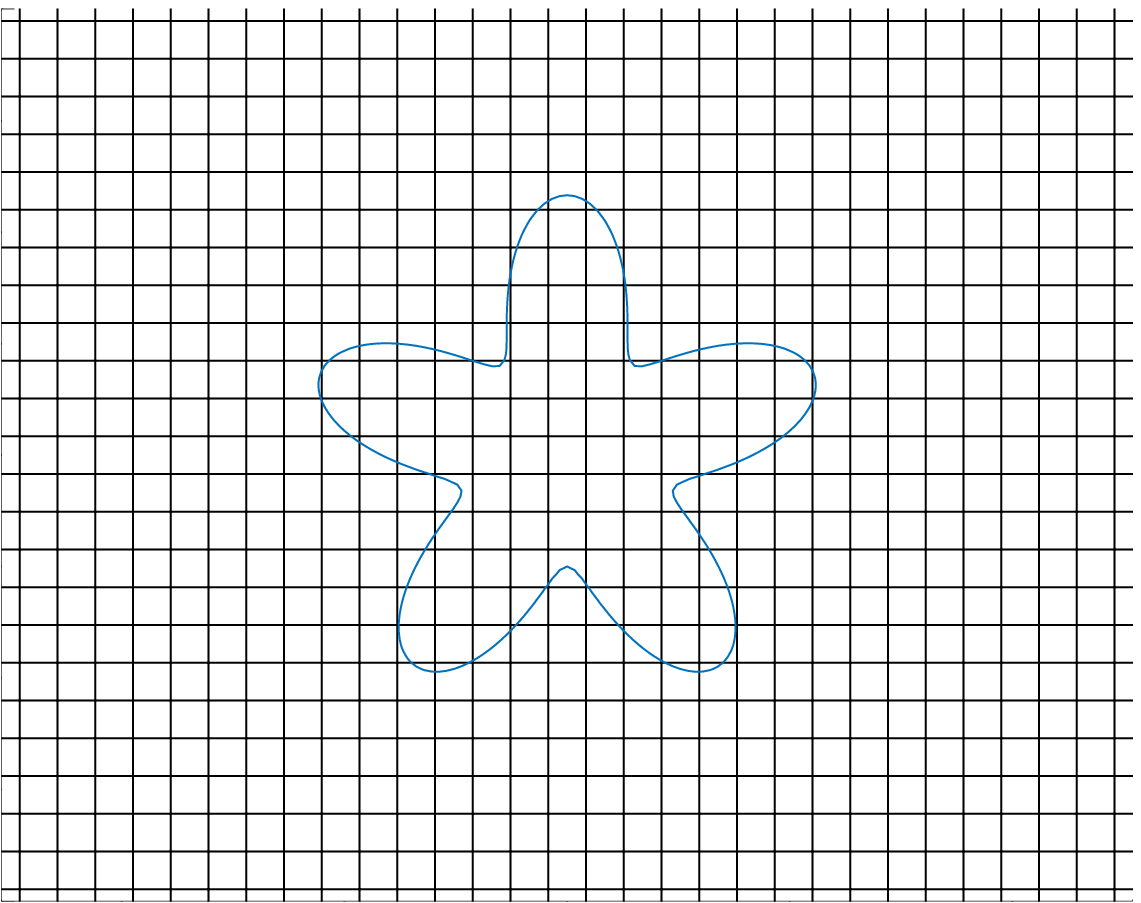}\\
{\scriptsize{}(a) Initial structured mesh and curved boundary $\Gamma$}
\par\end{center}%
\end{minipage}%
\begin{minipage}[t]{0.33\textwidth}%
\begin{center}
\includegraphics[width=0.9\textwidth]{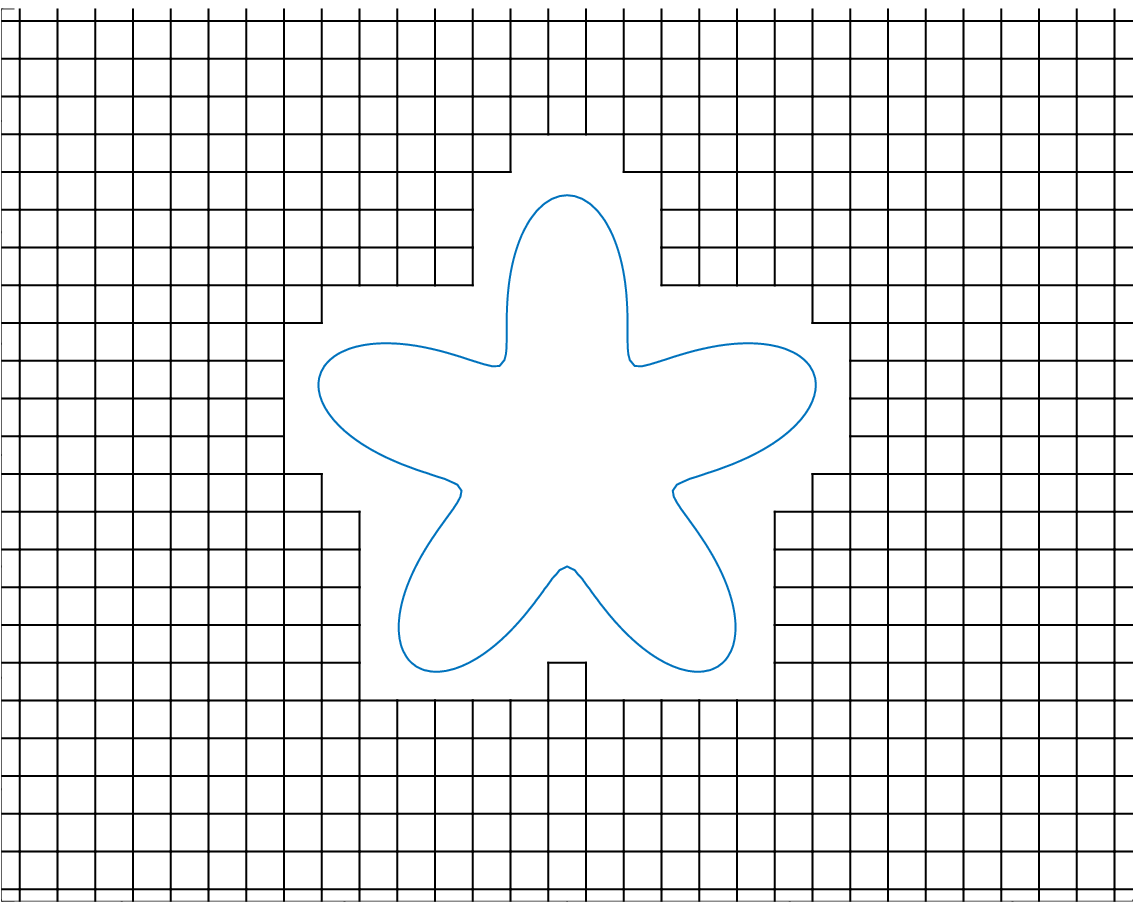}\\
{\scriptsize{}(b) After removing elements outside or on $\Gamma$}
\par\end{center}%
\end{minipage}%
\begin{minipage}[t]{0.33\textwidth}%
\begin{center}
\includegraphics[width=0.9\textwidth]{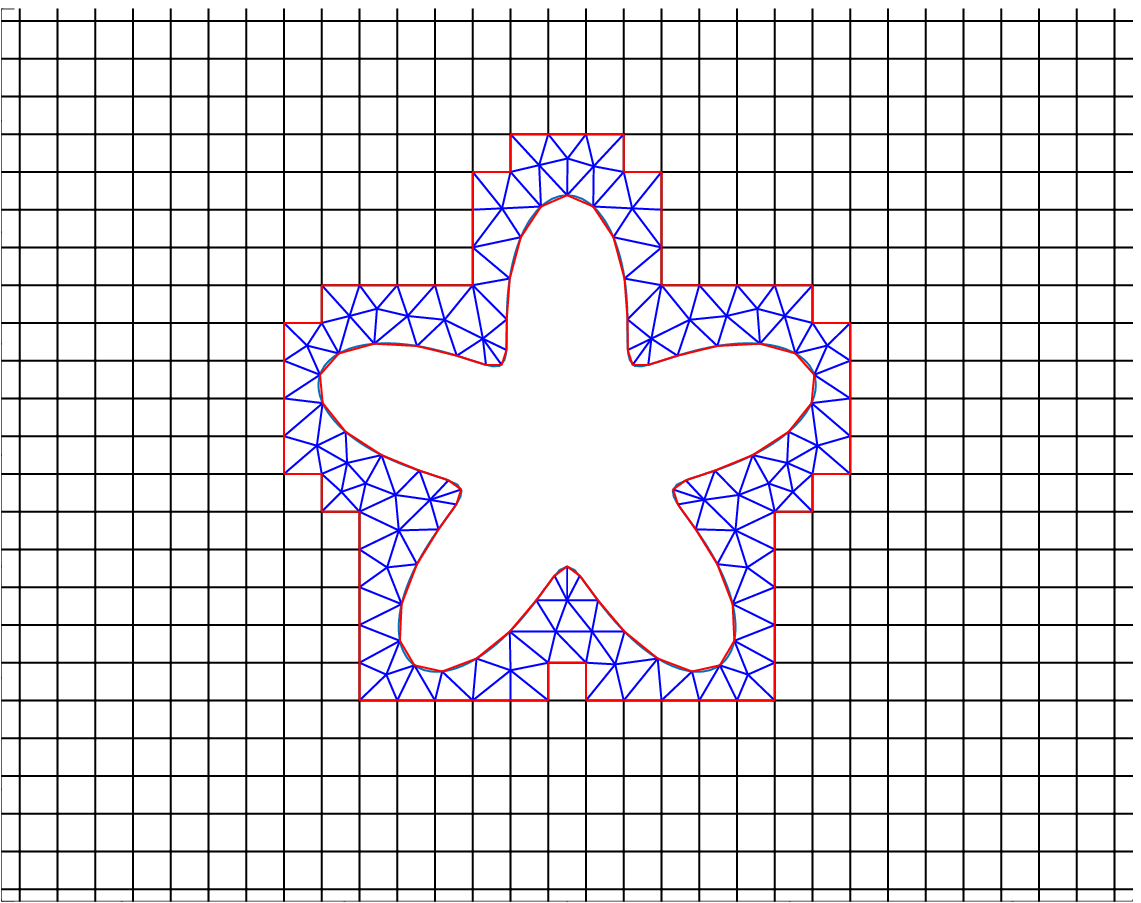}\\
{\scriptsize{}(c) After generating conformal triangles}
\par\end{center}%
\end{minipage}

\caption{\label{fig:initialmesh} Steps to generate initial linear mesh.}
\end{figure*}

To recover the high accuracy of SEM, we generate mixed-element meshes
using tensor-product elements in the interior of the geometric domain
as much as possible with triangular elements near curved boundaries.
To this end, we first generate a uniform quadrilateral mesh to cover
the bounding box of the domain. Next, we tag the nodes outside the
computational domain or within some fixed distance of a curved boundary,
where the distance is based on the average edge length of the initial
mesh. We then remove these nodes along with any quadrilateral containing
them. This process leaves a collection of connected edges that are
inside the geometric domain. We then mesh the gap between these edges
and the curved boundary using an off-the-shelf mesh generator. In
this work, we used MATLAB's built-in PDE Toolbox \cite{PDEToolbox},
of which the mesh generator built on top of CGAL \cite{cgal:eb-23a}.
This mesh generator allows nodes to be placed exactly on the curved
boundary. We apply $h$-GR adaptivity outlined in section~\ref{subsec:Applying-curvature-based-adaptiv}
during this phase. The nodes from the original quadrilateral mesh
and the new triangle mesh are matched so that those duplicate nodes
can be removed. The first couple of steps of the mesh generation are
shown in Figure~\ref{fig:initialmesh}. Next, the element connectivity
table is built using the half-facet array algorithm outlined in \cite{AHF}.
The resulting linear mesh might still have triangles with more than
one facet on the curved boundary. These elements could be problematic
as placing high-order edge nodes can cause the two facets to ``flatten
out'' and create a negative Jacobian determinant at the corner node
where the two curved edges intersect. We use a simple edge flip to
prevent this case and ensure that all triangles have only one boundary
facet. This process gives a good-quality linear finite element mesh.
However, if stronger mesh quality constraints are required, we can
also further improve it by using mesh-smoothing techniques, such as
that in \cite{jiao2011simple}.

\subsection{Applying $h$-GR\label{subsec:Applying-curvature-based-adaptiv}}

\begin{figure*}
\begin{minipage}[t]{0.5\textwidth}%
\begin{center}
\includegraphics[width=0.6\textwidth]{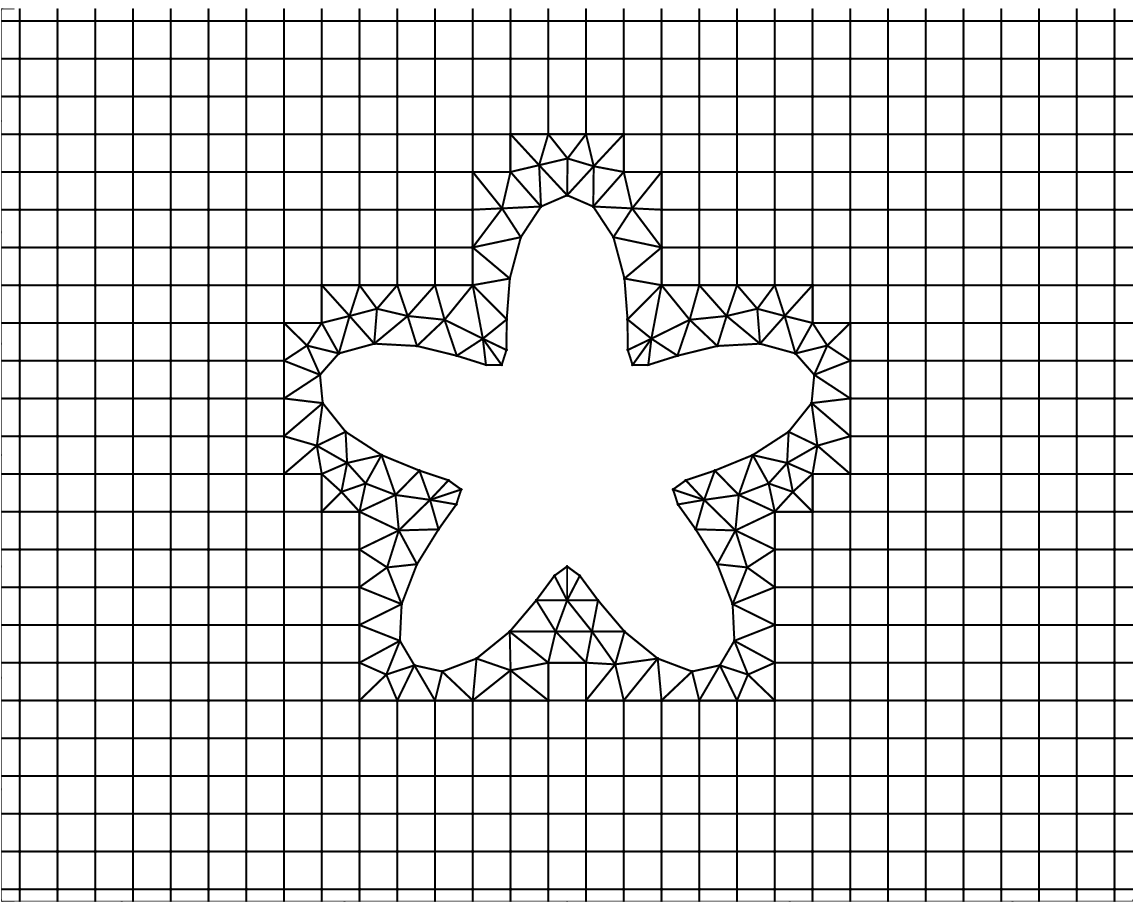}\\
{\scriptsize{}(a) Mesh without $h$-GR}
\par\end{center}%
\end{minipage}%
\begin{minipage}[t]{0.5\textwidth}%
\begin{center}
\includegraphics[width=0.6\textwidth]{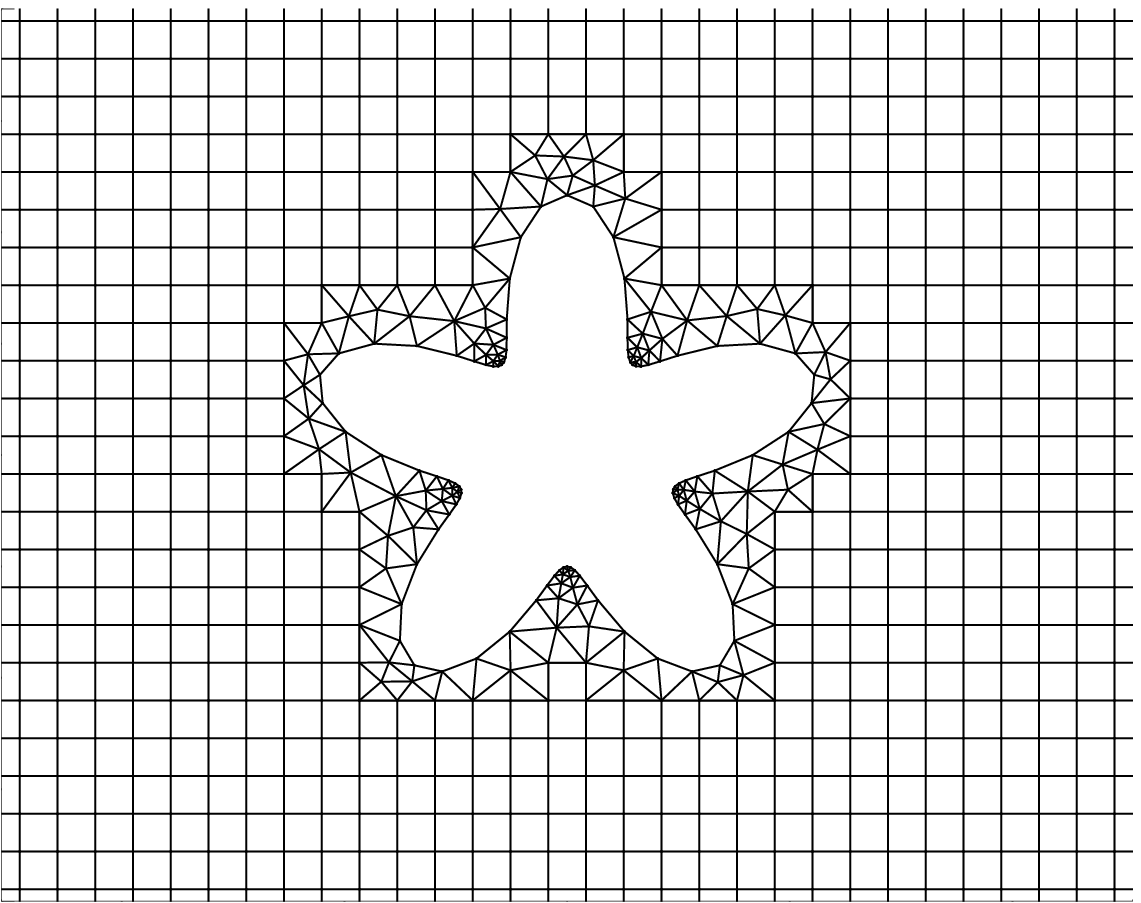}\\
{\scriptsize{}(b) Mesh with $h$-GR}
\par\end{center}%
\end{minipage}\caption{\label{fig:CBRmeshes} Example meshes around a flower hole with and
without curvature-based refinement.}
\end{figure*}

To apply $h$-GR, we seek to generate a local edge length $h$ at
each node on the boundary, which is linearly related to the curvature
of the boundary. This step is relatively simple and adds negligible
computational cost to the mesh generation. We evaluate the curvature
for each node on the curved boundary using either the explicit parameterization
or some form of geometric reconstruction. Next, we define a target
angle $\theta_{max}$ such that no two adjacent nodes have a difference
in the normal direction that is greater than $\theta_{max}$. To this
end, we make the target edge length the arc length of the circle defined
by the radius of curvature and $\theta_{max}$. The edge length at
node $i$ is based on the radius of curvature and can be defined as
$h_{i}=\theta_{max}R_{i}=\theta_{max}/K_{i}$ where $R_{i}$ is the
radius of curvature and $K_{i}$ is the unsigned curvature for node
$i$. Since $h_{i}$ defined in this way can be arbitrarily large
or small, we limit the minimum and maximum edge lengths using parameters
$h_{min}$ and $h_{max}$, respectively. The overall formula for computing
the target edge length at node $i$ is defined as
\[
h_{i}=\mathrm{min}(\mathrm{max}(\theta_{max}/K_{i},h_{min}),h_{max}).
\]
These edge lengths are supplied to the mesh generator as target edge
lengths for generating the mesh near the boundary. A mesh with and
without $h$-GR adaptivity is displayed in Figure~\ref{fig:CBRmeshes}.

\subsection{Placing High-Order Nodes\label{subsec:Placing-High-Order-Nodes}}

After obtaining the initial linear mesh, we insert mid-edge and mid-face
nodes to generate high-order spectral elements. For the rectangular
elements in the interior of the domain, we insert Gauss-Lobatto points
to generate tensor-product spectral elements. For the near-boundary
elements, we also use Gauss-Lobatto nodes for the facets (i.e., edges)
shared between tensor-product and non-tensor-product elements to ensure
continuity. The placements of mid-face nodes of triangles are more
flexible; we place them at quadrature points that maximize the degree
of the quadrature rule as Gauss-Lobatto points do for tensor-product
elements. If an element has a curved edge, however, the placement
of the nodes needs to satisfy the so-called Ciarlet-Raviart condition
\cite{ciarlet1972combined} to preserve the convergence rate. In a
nutshell, the Ciarlet-Raviart condition requires the Jacobian determinant
of a degree-$p$ element to have bounded derivatives up to $(p+1)$st
order \cite{li2019compact}. To satisfy the condition, we use an iterative
process similar to the one presented in \cite{li2019compact}. In
particular, to generate high-order curved elements with degree-$p$
polynomial basis functions, starting from $q=2$, we incrementally
generate a degree-$q$ element by first interpolating the nodal positions
using a degree-$(q-1$) element and then project the nodes on curved
edges onto the exact geometry. We repeat the process until $q=p+2$,
then interpolate the nodal positions of the degree-$p$ element from
the degree-$(p+2)$ element. In this fashion, we found that the Jacobian
determinant is sufficiently smooth to reduce the pollution errors
of the spectral elements. Note that this process is needed only for
elements directly incident on curved boundaries, so its computational
cost is less complex than the rest of the mesh-generation process.

\subsection{Extraction of near boundary elements for AES-FEM\label{subsec:near-boundary-aes-fem}}

To utilize AES-FEM near the boundary, we need a linear mesh. To this
end, we construct a new, smaller, linear mesh from pieces of the high-order
mesh. We first decompose the triangular part of the mesh and include
all newly generated linear triangles in the AES-FEM mesh. Choosing
just the decomposed triangles from the original mesh will not be sufficient
to improve accuracy. We must extend into the quadrilateral section
of the mesh by a number of layers. A layer is defined as all quadrilaterals
with two or more nodes as part of the AES-FEM mesh. 

\begin{figure}
\begin{minipage}[t]{0.45\textwidth}%
\begin{center}
\includegraphics[scale=0.45]{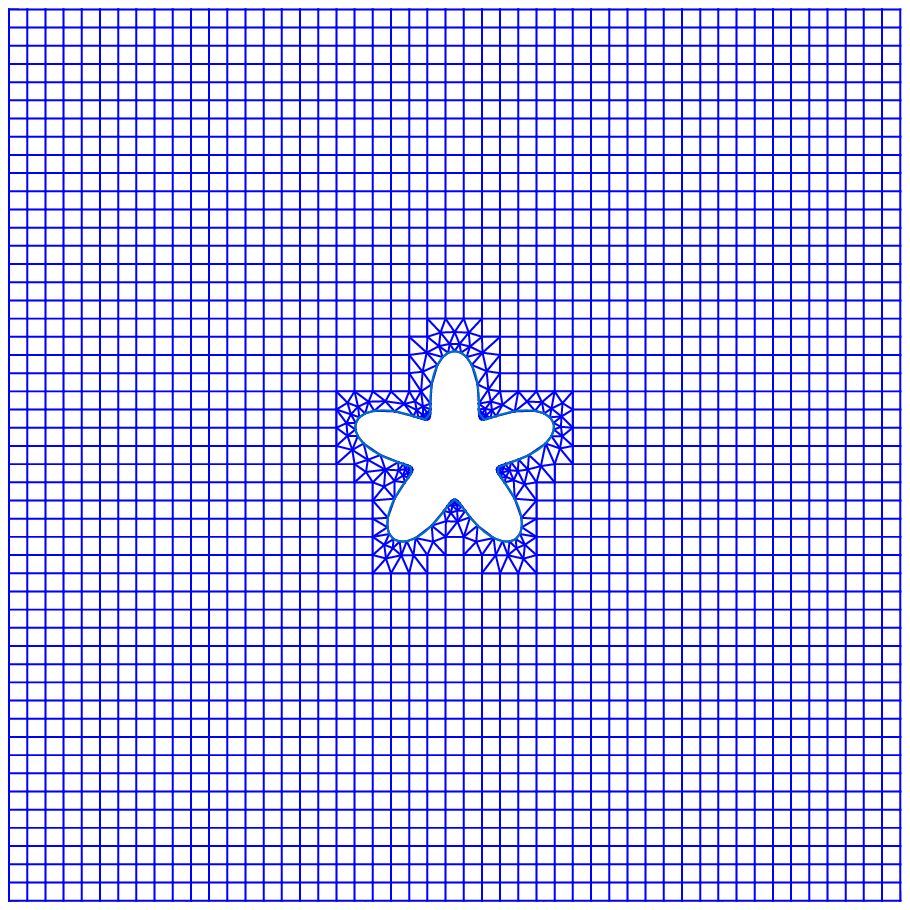}
\par\end{center}
\begin{center}
{\scriptsize{}(a) Quadratic mesh before splitting elements near the
boundary}{\scriptsize\par}
\par\end{center}%
\end{minipage}%
\begin{minipage}[t]{0.45\textwidth}%
\begin{center}
\includegraphics[scale=0.45]{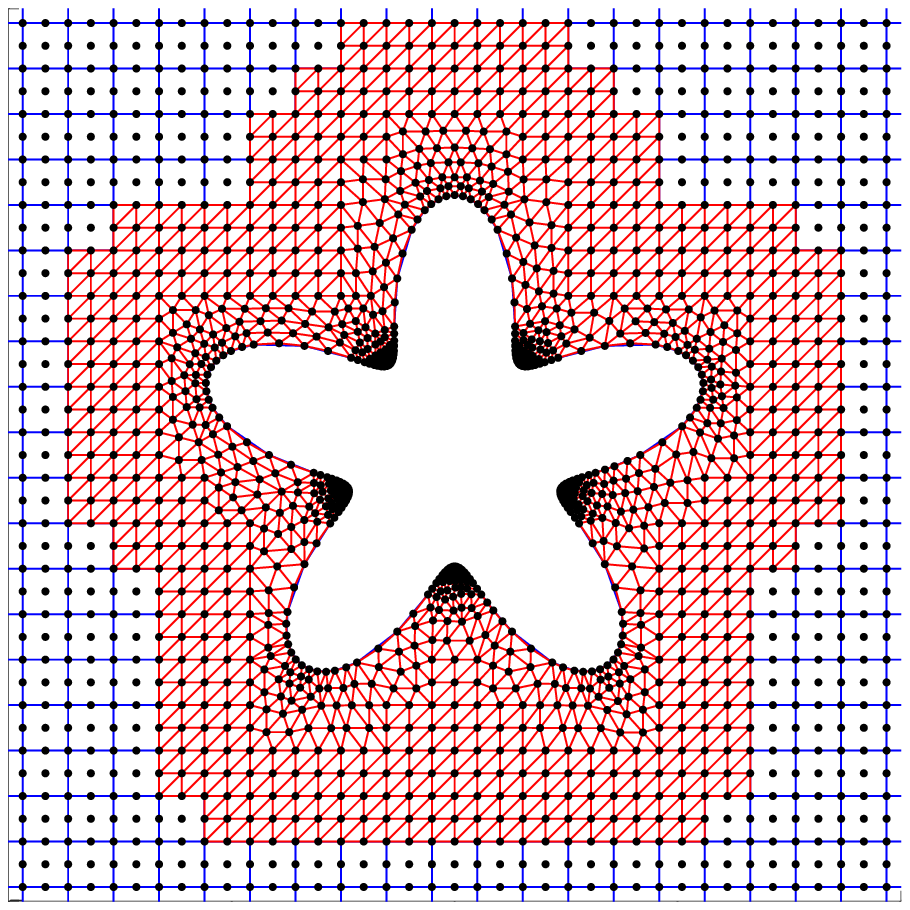}
\par\end{center}
\begin{center}
{\scriptsize{}(b) AES-FEM boundary mesh after extraction and splitting}{\scriptsize\par}
\par\end{center}%
\end{minipage}

\caption{\label{fig:Aes-extract-mesh} Example meshes before and after splitting
near the boundary. In (b), the red linear triangles constitute the
AES-FEM mesh, obtained by splitting quadratic elements.}

\end{figure}

For this reason, we keep a boolean array of all nodes currently in
the AES-FEM mesh. Care must be taken as too few elements will result
in only minor improvement by AES-FEM, while too many elements will
waste computing power (although it will be stable). In practice, two
layers typically suffice, and one layer may also suffice for quadratic
meshes of simple domains such as the elliptical hole domain. After
all the quadrilaterals have been identified, they are decomposed into
linear elements and then split into triangles and added to the AES-FEM
mesh. Figure~\ref{fig:Aes-extract-mesh} shows an example of an extracted
boundary mesh. To perform post-processing, we then use AES-FEM to
solve the PDE by using the SEM solutions as Dirichlet boundary conditions
for the artificial boundary.

\subsection{Extension to Three Dimensions}

While the process outlined in the paper is for 2D, we can extend this
mesh generation method to 3D by addressing some additional challenges.
First, unstructured mesh generation in 3D uses tetrahedra as the primary
element shape. Superconvergent spectral elements, however, require
the use of tensor-product elements, which, in 3D, have quadrilateral
faces. To construct a conformal mesh, we must resolve the two boundaries
by introducing pyramid-shaped elements, which cover the quadrilateral
face boundary of the hexahedral mesh and the triangular boundary of
the tetrahedral mesh. Meanwhile, the curved surface can be meshed
using any standard surface mesh technique. The tetrahedral mesh generator
meshes the space between the two boundaries. The fine-tuning for this
part is determining the minimum distance from the pyramids' boundary
to ensure the stability of the tetrahedral mesh generation. Curvature
in this context can be replaced by the maximum curvature \cite{Carmo1976}.
More specifically, given a smooth parametric surface $S(u,v)\in\mathbb{R}^{2}\mapsto\mathbb{R}^{3}$,
let $\boldsymbol{J}$ be the Jacobian matrix of $S$. We define the
first fundamental matrix of $S$ as $\boldsymbol{G}=\boldsymbol{J}^{T}\boldsymbol{J}$.
Since the vectors of the Jacobian matrix form a basis for the tangent
space of the surface, we can define the normal direction of the surface
as 
\[
\hat{\boldsymbol{n}}=\frac{S_{u}\times S_{v}}{\parallel S_{u}\times S_{v}\parallel}.
\]
The second fundamental matrix is then defined as 
\[
\boldsymbol{B}=\left[\begin{array}{cc}
\hat{\boldsymbol{n}}^{T}S_{uu} & \hat{\boldsymbol{n}}^{T}S_{uv}\\
\hat{\boldsymbol{n}}^{T}S_{vu} & \hat{\boldsymbol{n}}^{T}S_{vv}
\end{array}\right].
\]
The shape operator is then $\boldsymbol{W}=\boldsymbol{G}^{-1}\boldsymbol{B}$,
of which the larger absolute value of its eigenvalues is the maximum
curvature; the detailed computation can be found \cite{Wang2009}.
We then use its inverse, namely the minimal radius of curvature, as
the criteria for controlling the local edge lengths similar to the
process outlined in section~\ref{subsec:Applying-curvature-based-adaptiv}.
As for the rest of the 2D algorithm, section~\ref{subsec:Placing-High-Order-Nodes}
generalizes to 3D directly. We defer the implementation of 3D to future
work, and we will focus on verifying the error analysis and the feasibility
of the overall workflow in 2D.

\section{Numerical Results\label{sec:Results}}

\begin{figure*}
\begin{centering}
\begin{minipage}[t]{0.5\textwidth}%
\begin{center}
\includegraphics[width=0.8\columnwidth]{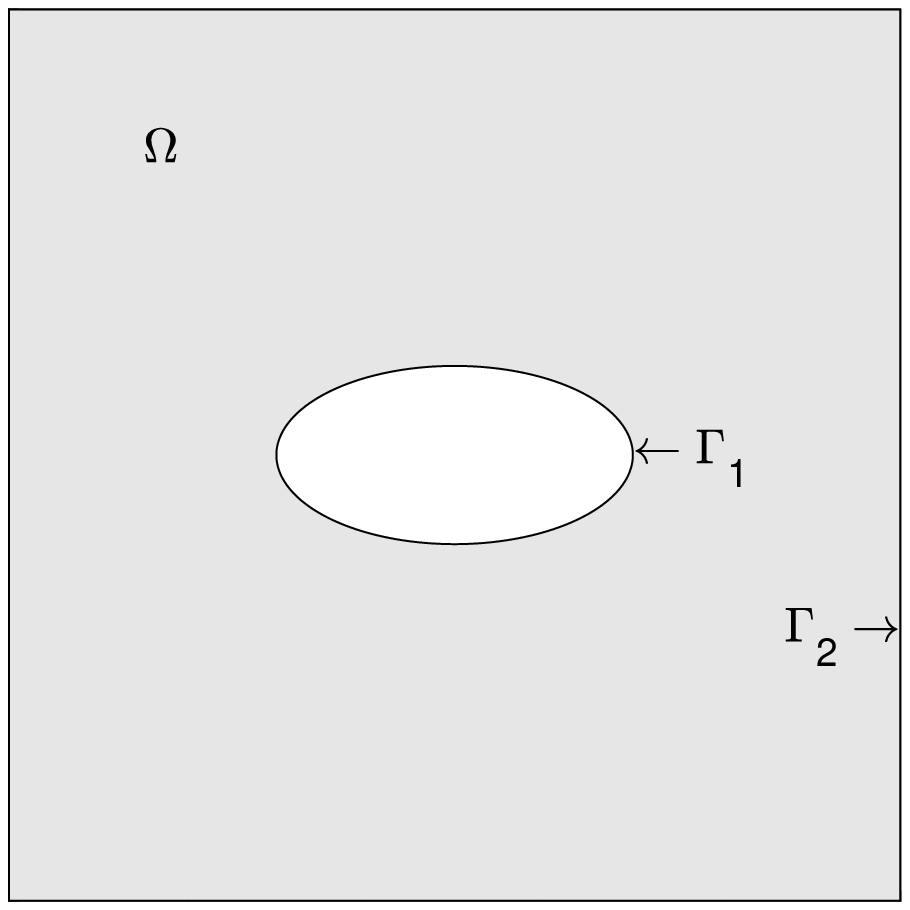}\\
{\scriptsize{}(a) Unit square with elliptical hole}
\par\end{center}%
\end{minipage}%
\begin{minipage}[t]{0.5\textwidth}%
\begin{center}
\includegraphics[width=0.8\columnwidth]{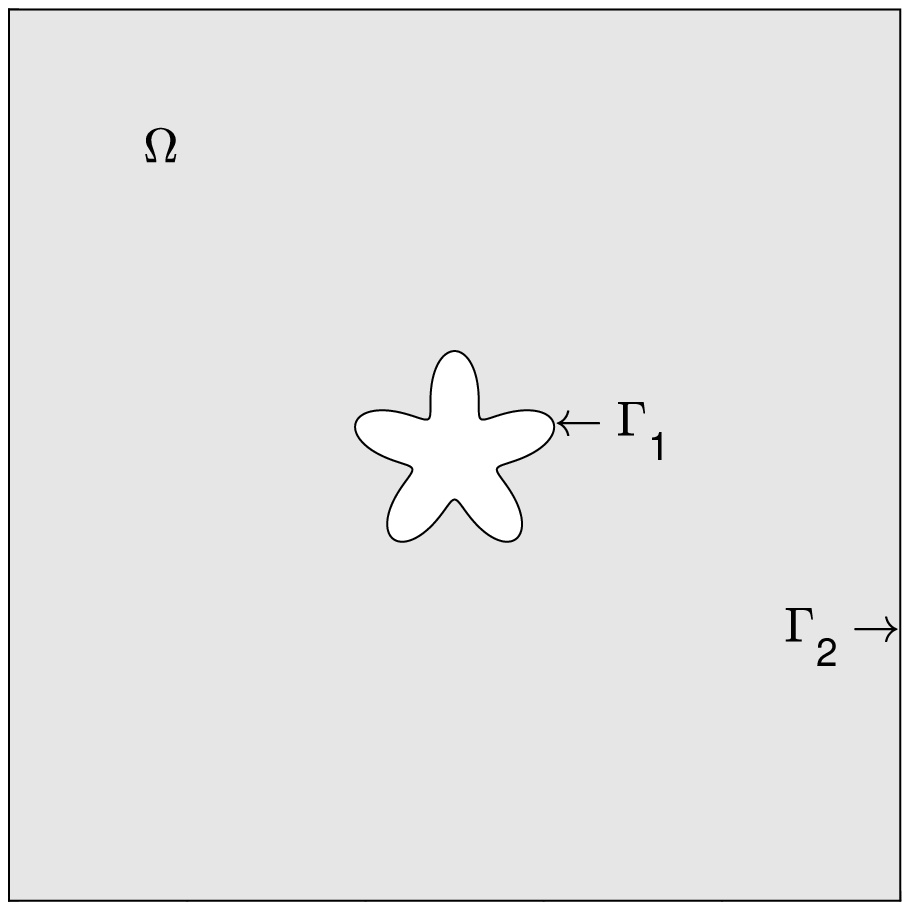}\\
{\scriptsize{}(b) Unit square with five-petal flower hole}
\par\end{center}%
\end{minipage}
\par\end{centering}
\caption{\label{fig:domain_figures}Domains used in tests. Neumann or Dirichlet
boundary conditions are applied to interior boundary $\Gamma_{1}$;
Dirichlet boundary conditions are applied to exterior boundary $\Gamma_{2}$.}
\end{figure*}

In this section, we present numerical results with our proposed method,
primarily as a proof of concept in improving the numerical accuracy
of SEM over curved domains. To this end, we solve the convection-diffusion
equation
\begin{align}
-\Delta u+\boldsymbol{v}\cdot\boldsymbol{\nabla}u=f & \quad\text{on }\Omega,\label{eq:conv-diff}\\
u=u_{D} & \quad\text{on }\Gamma_{D},\label{eq:conv-diff-Dirichlet}\\
\partial_{\boldsymbol{n}}u=u_{N} & \quad\text{on }\Gamma_{N},\label{eq:conv-diff-Neumann}
\end{align}
where $f$ denotes a source term, and $u_{D}$ and $u_{N}$ denote
the Dirichlet and Neumann boundary conditions, respectively. To assess
accuracy, we used the method of manufactured solutions in 2D and computed
$f$, $u_{D}$, and $u_{N}$ from the exact solution $u=\sin(10\pi x)\mathrm{cos}(10\pi y)+xy$
and $\boldsymbol{v}=[x,-y].$ We tested the method using two domains:
a unit square centered at $(0.5,0.5)$ with an elliptical hole and
a five-petal flower hole, respectively, as illustrated in Figure~\ref{fig:domain_figures}.
For the elliptical hole, the semi-major and semi-minor axes are 0.2
and 0.1, respectively. The five-petal flower is a variant of the flower-shaped
domain in \cite{cheung2019optimally}, with the parametric equations
\begin{align*}
x= & (0.25+0.1\sin(5\theta))\cos(\theta)/3+0.5\\
y= & (0.25+0.1\sin(5\theta))\sin(\theta)/3+0.5
\end{align*}
for $0\leq\theta<2\pi$. This domain is particularly challenging due
to its high curvature and mixed convex-and-concave regions.

To demonstrate superconvergence, we measure the errors in $\ell_{2}$
norm instead of $L_{2}$ norm because the latter is dominated by the
interpolation error. More specifically, let $\boldsymbol{u}$ and
$\hat{\boldsymbol{u}}$ denote the vectors containing the exact and
numerical solutions at the nodes; the $\ell_{2}$-norm error is then
computed as
\[
\mathrm{error}=\Vert\boldsymbol{u}-\hat{\boldsymbol{u}}\Vert_{2}=\sqrt{\sum_{i}(u_{i}-\hat{u}_{i})^{2}}.
\]
Since the elements near the boundaries are irregular, we compute the
convergence rate based on the number of degrees of freedom (dof),
i.e., 
\[
\mathrm{convergence\:rate}=-\log\left(\frac{\mathrm{error\,of\,fine\,mesh}}{\mathrm{error\,of\,coarse\,mesh}}\right)\biggl/\log\left(\sqrt[d]{\frac{\mathrm{dof\,in\,fine\,mesh}}{\mathrm{dof\,in\,coarse\,mesh}}}\right),
\]
where $d=2$ is the topological dimension of the mesh. This formula
is consistent with the edge-length-based convergence rate for structured
meshes for sufficiently fine meshes. 

\subsection{Effect of $h$-GR on Curved Neumann Boundaries\label{subsec:h-refinement}}

\begin{figure*}
\begin{minipage}[t]{0.45\textwidth}%
\begin{center}
\includegraphics[width=1\textwidth]{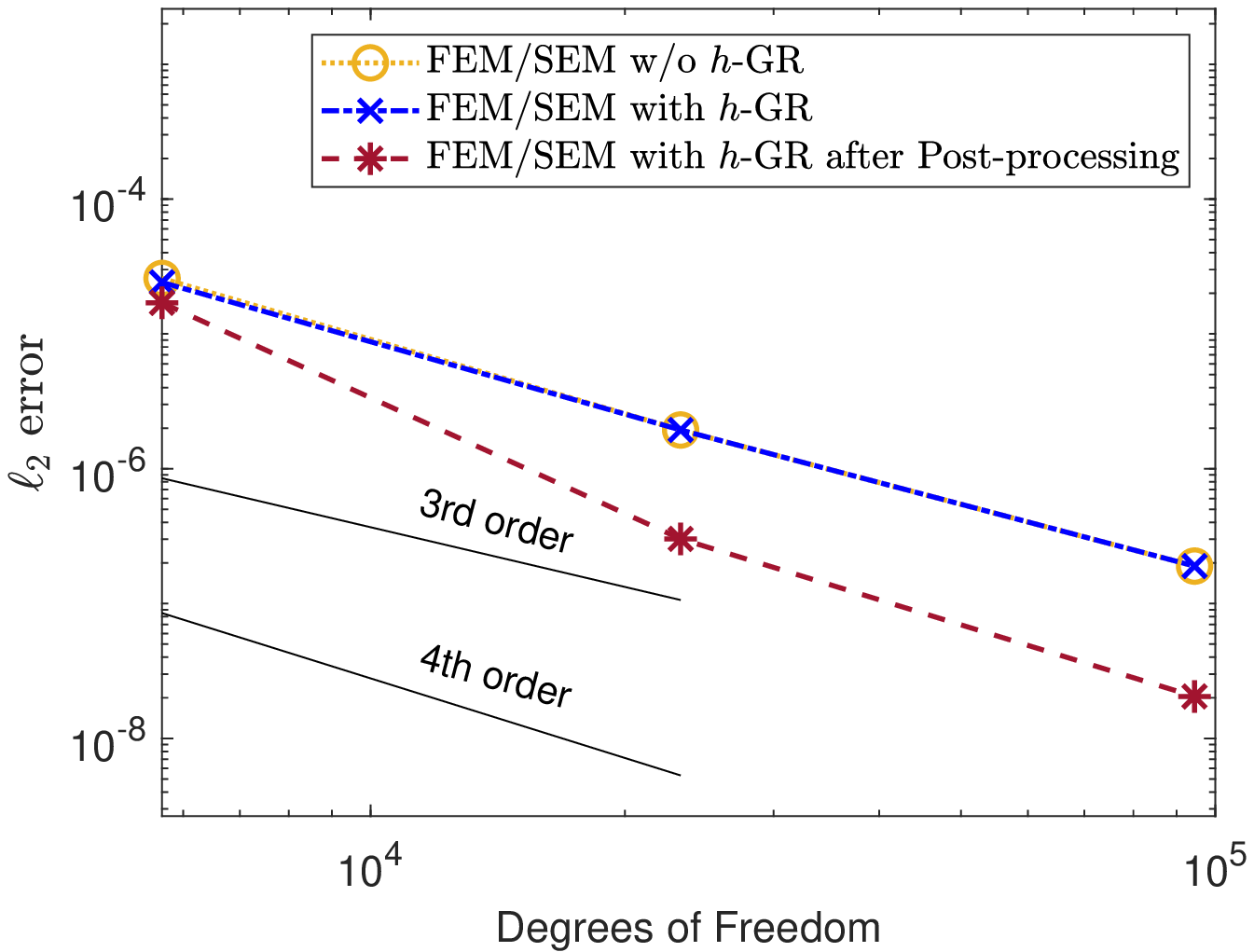}\\
{\scriptsize{}(a) Quadratic elliptical hole}
\par\end{center}%
\end{minipage} %
\begin{minipage}[t]{0.45\textwidth}%
\begin{center}
\includegraphics[width=1\textwidth]{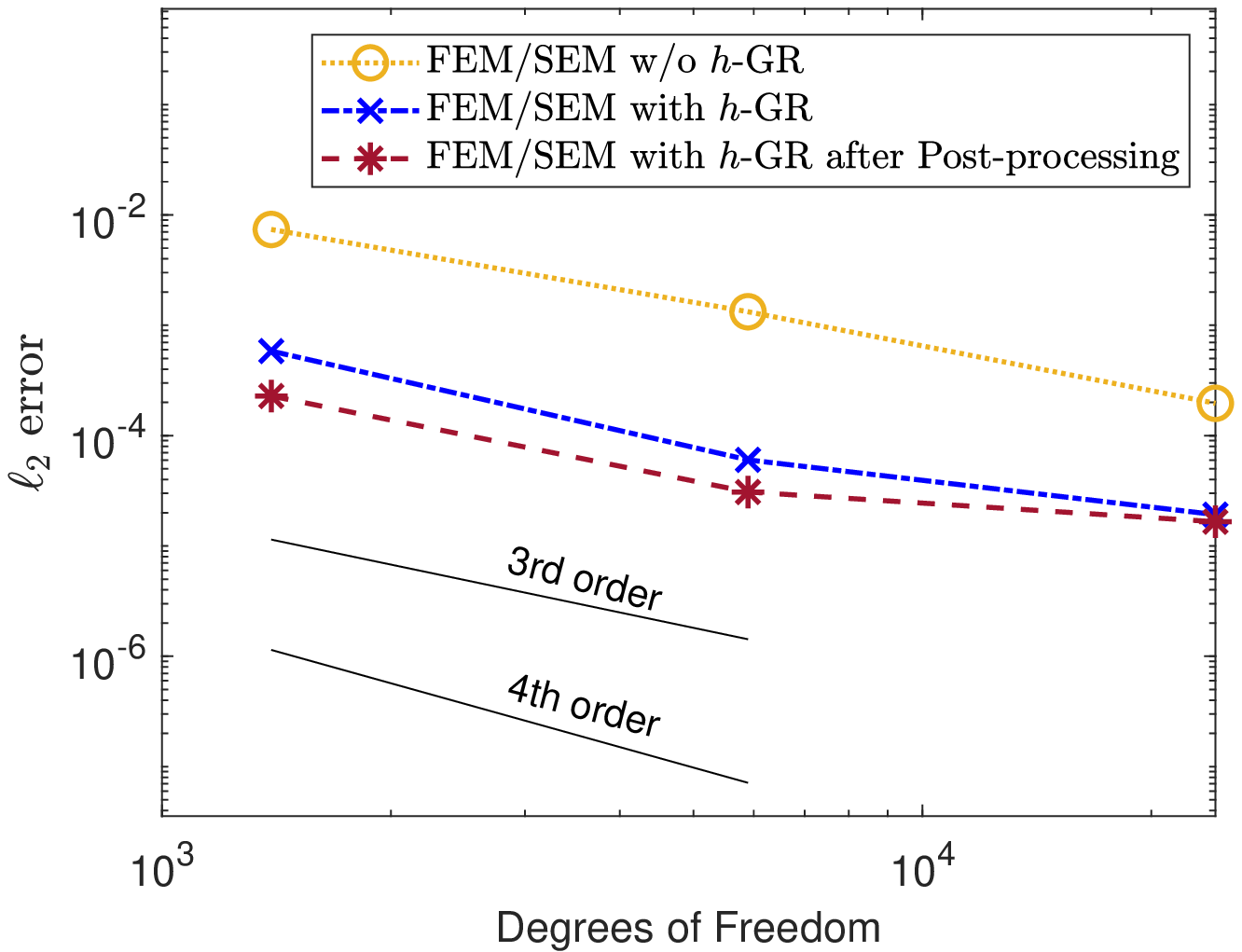}\\
{\scriptsize{}(b) Quadratic flower hole}
\par\end{center}%
\end{minipage}

\begin{minipage}[t]{0.45\textwidth}%
\begin{center}
\includegraphics[width=1\textwidth]{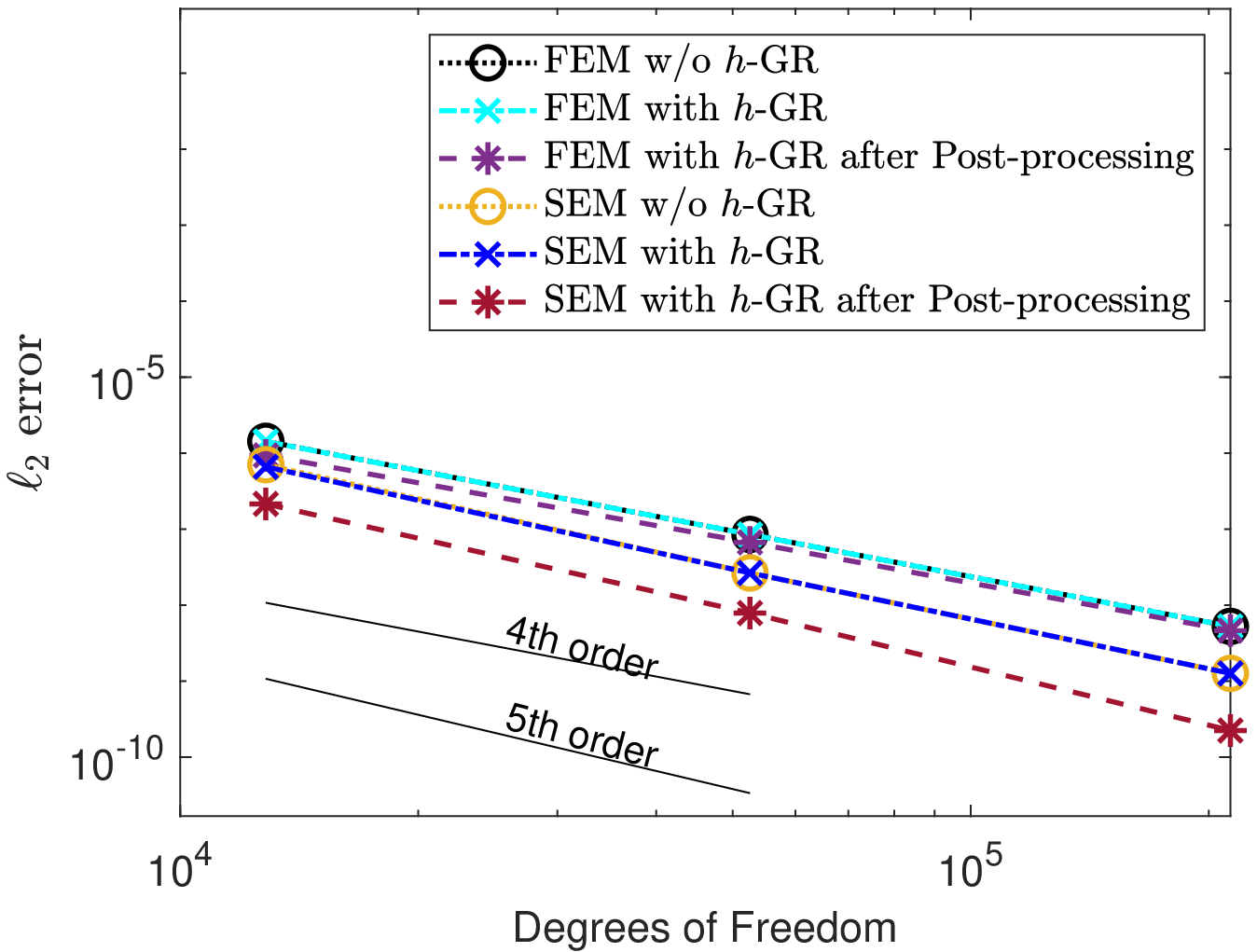}\\
{\scriptsize{}(c) Cubic elliptical hole}
\par\end{center}%
\end{minipage} %
\begin{minipage}[t]{0.45\textwidth}%
\begin{center}
\includegraphics[width=1\textwidth]{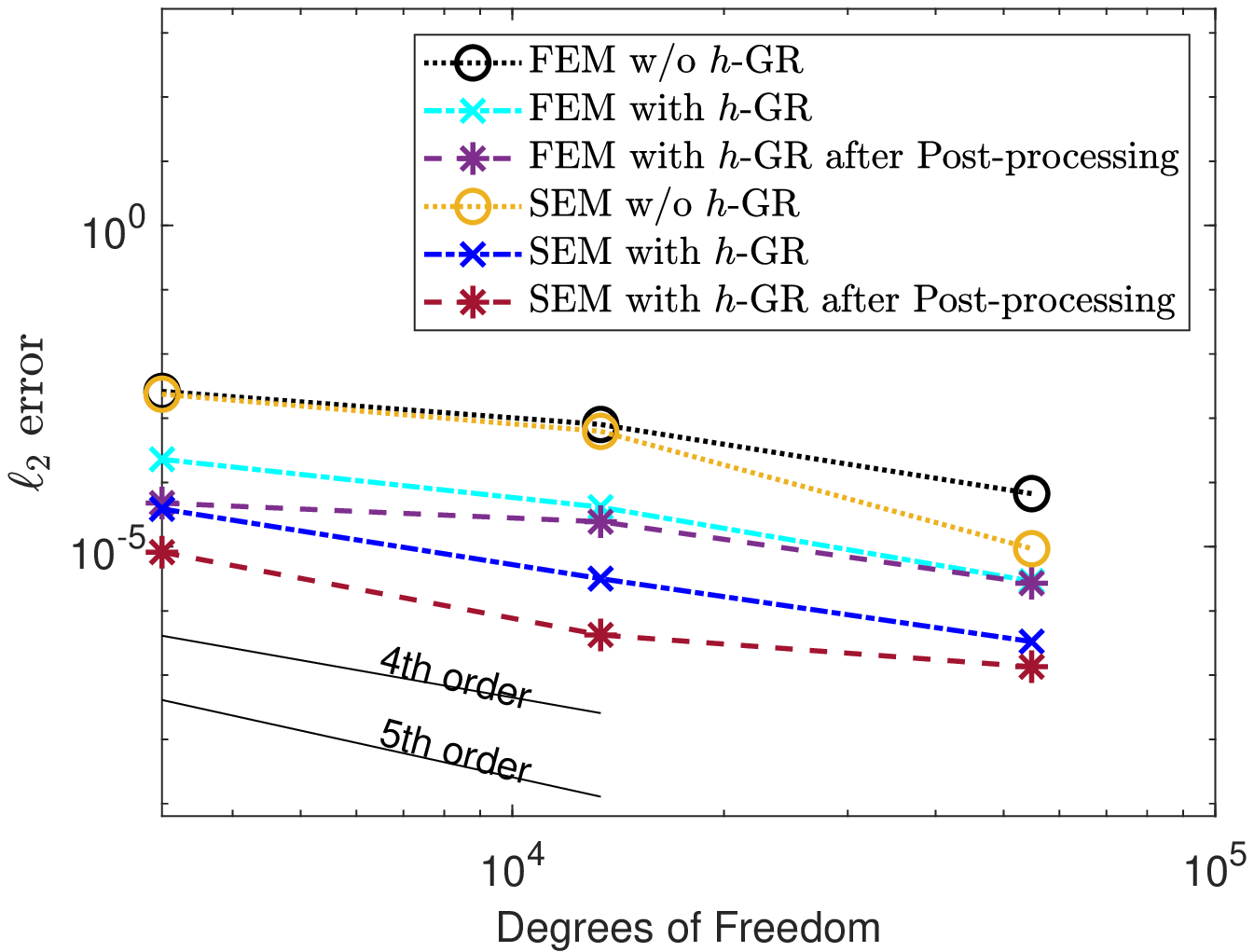}\\
{\scriptsize{}(d) Cubic flower hole}
\par\end{center}%
\end{minipage}

\begin{minipage}[t]{0.45\textwidth}%
\begin{center}
\includegraphics[width=1\textwidth]{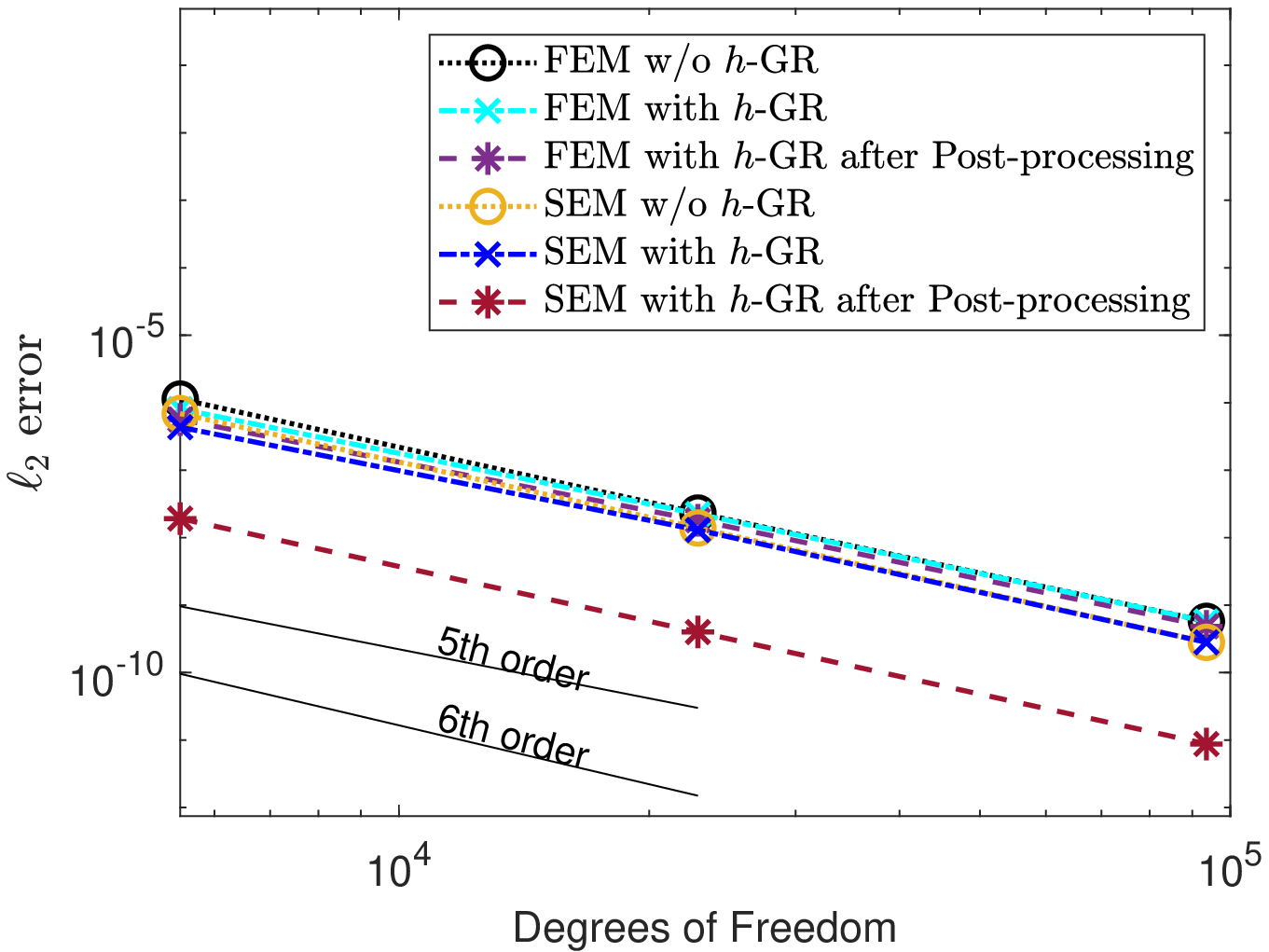}\\
{\scriptsize{}(e) Quartic elliptical hole}
\par\end{center}%
\end{minipage} %
\begin{minipage}[t]{0.45\textwidth}%
\begin{center}
\includegraphics[width=1\textwidth]{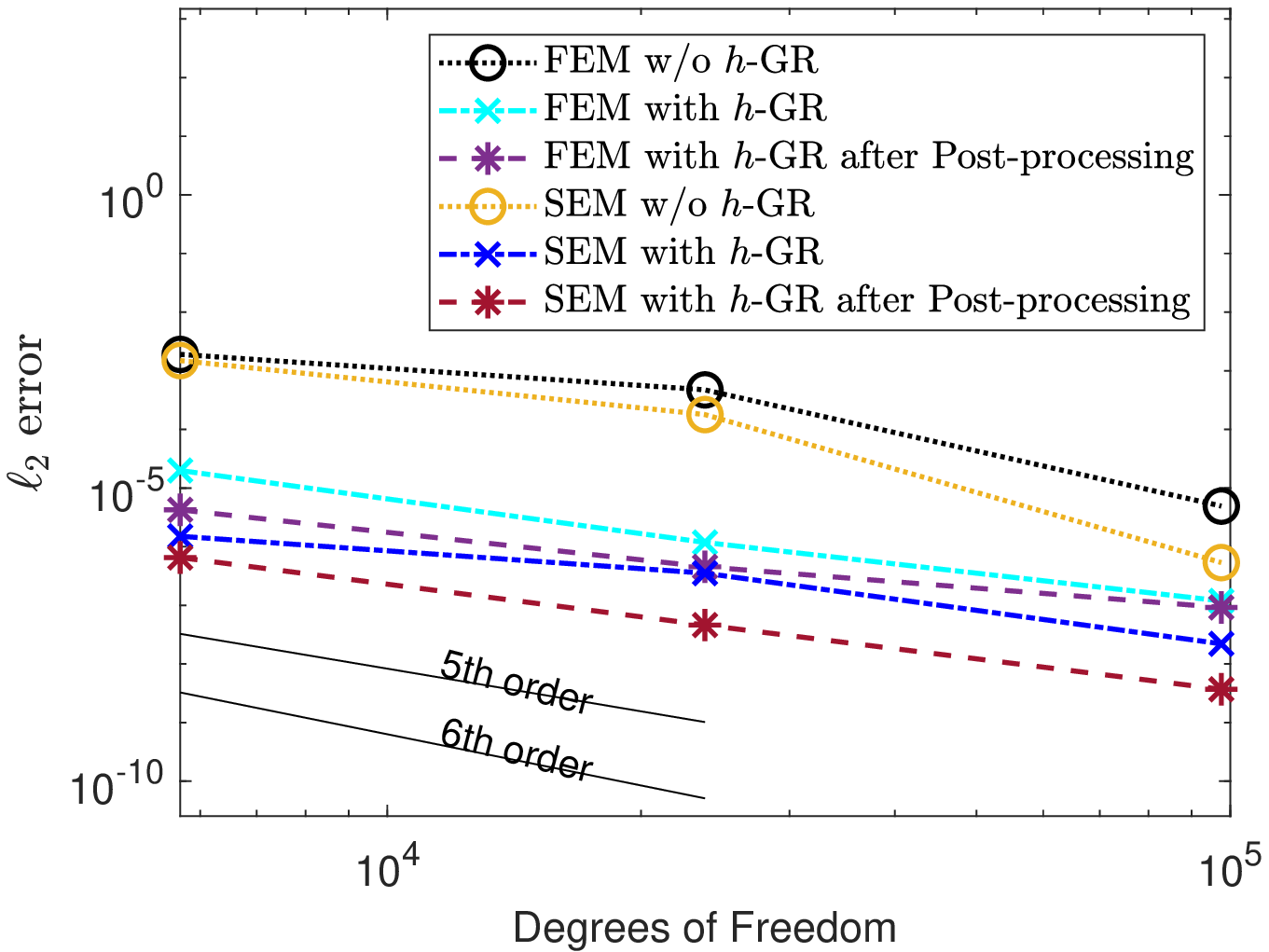}\\
{\scriptsize{}(f) Quartic flower hole}
\par\end{center}%
\end{minipage}

\caption{\label{fig:CBR_ell_iso}Comparison of nodal errors of spectral elements
in the interior of the flower-hole domain with quadratic, cubic, and
quartic SEM for the convection-diffusion equation with and without
$h$-GR, as well as the error after post-processing.}
\end{figure*}

We first assess the effect of $h$-GR on curved Neumann boundaries.
To this end, we solved the convection-diffusion equation (\ref{eq:conv-diff})
on the elliptical-hole and flower-hole meshes with Neumann boundary
conditions on the curved boundary $\Gamma_{1}$ of the flower hole
and Dirichlet boundary conditions on $\Gamma_{2}$, as illustrated
in Figure~\ref{fig:domain_figures}(a) and (b), respectively. We
solved the equation using isoparametric FEM and SEM with and without
$h$-GR; we studied the effect of superparametric elements in section~\ref{subsec:hp-refinement}.
Figure~\ref{fig:CBR_ell_iso} shows the $\ell_{2}$-norm errors of
the nodal solutions with and without $h$-GR for quadratic, cubic,
and quartic FEM and SEM. For completeness, Figure~\ref{fig:CBR_ell_iso}
also shows the results after post-processing the near-boundary elements.
Note that quadratic FEM and SEM are identical, so there is only one
set of results for both. 

From the results, we first observe that cubic and quartic SEM significantly
outperformed their FEM counterparts in all cases. Hence, we will consider
only SEM in our later discussions. We also observe that for the flower-hole
mesh, $h$-GR alone improved the overall accuracy by about one, two,
and three orders of magnitude for quadratic, cubic, and quartic SEM
on the coarser and intermediate meshes. This drastic improvement is
because $h$-GR improves the accuracy of the normal directions of
high-curvature regions. The normals are implicitly used by Neumann
boundary conditions in (\ref{eq:BVP Neumann}). In contrast, for the
elliptical-hole mesh, $h$-GR alone did not improve accuracy much
because the curvatures are not high for the ellipse, and the normals
are reasonably accurate even without $h$-GR. Nevertheless, $h$-GR
still improved the accuracy of the spectral elements near boundaries,
even for the elliptical-hole domain. This improvement is evident because
the post-processed SEM solution was significantly better with $h$-GR
than without $h$-GR.

\begin{figure*}
\begin{centering}
\begin{minipage}[t]{0.45\textwidth}%
\begin{center}
\includegraphics[width=1\textwidth]{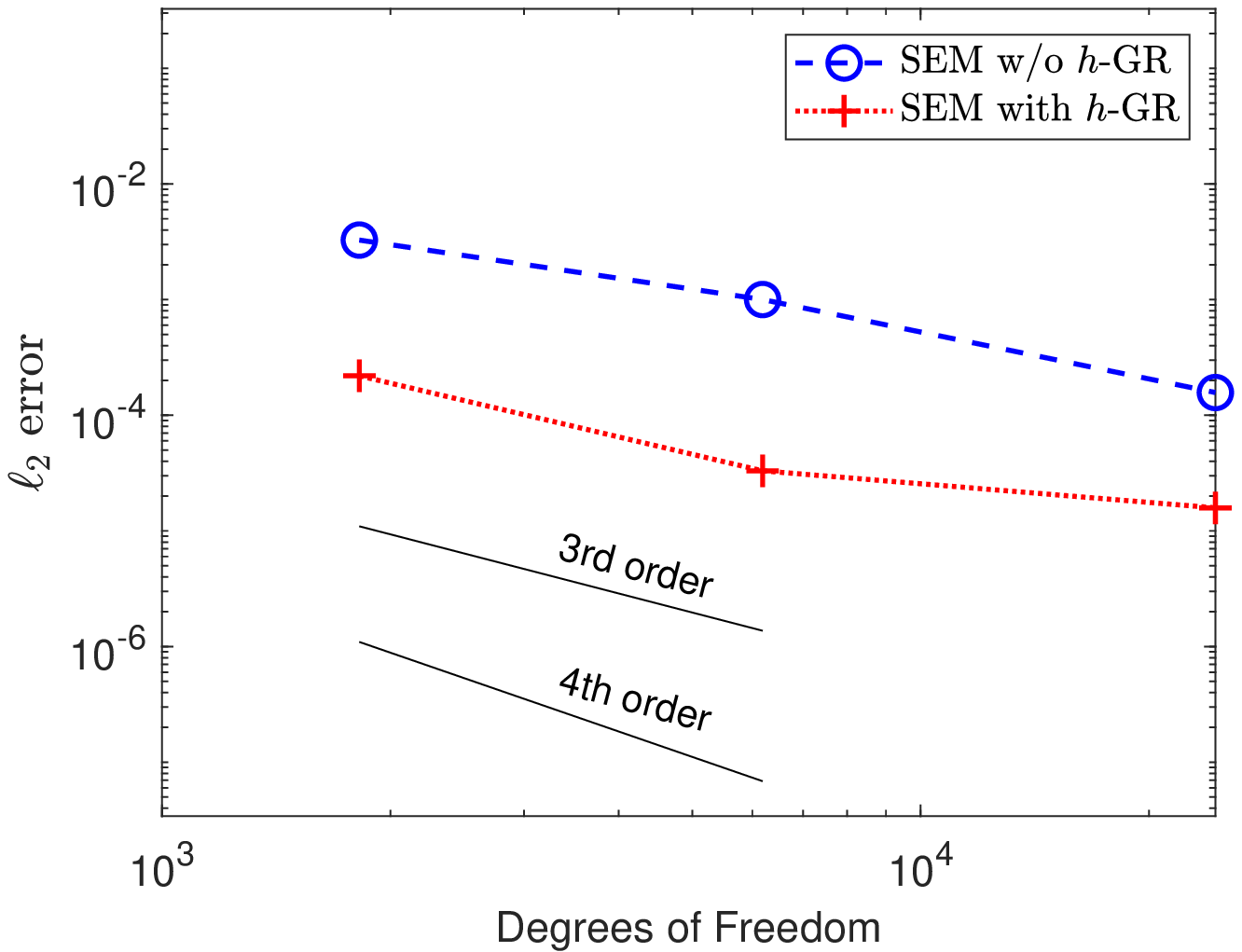}\\
{\scriptsize{}(a) Quadratic}
\par\end{center}%
\end{minipage}%
\begin{minipage}[t]{0.45\textwidth}%
\begin{center}
\includegraphics[width=1\textwidth]{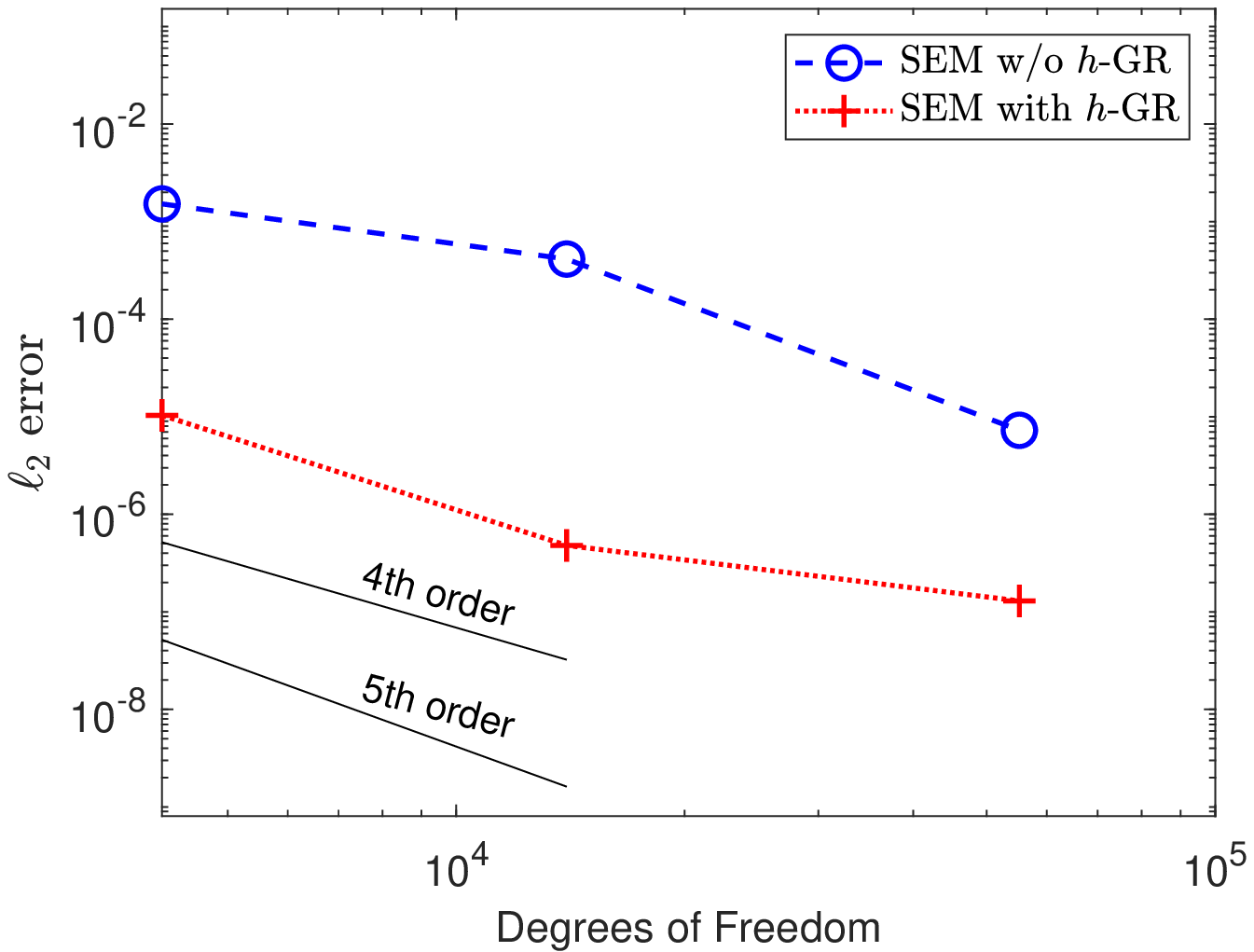}\\
{\scriptsize{}(b) Cubic}
\par\end{center}%
\end{minipage}
\par\end{centering}
\begin{centering}
\begin{minipage}[t]{0.45\textwidth}%
\begin{center}
\includegraphics[width=1\textwidth]{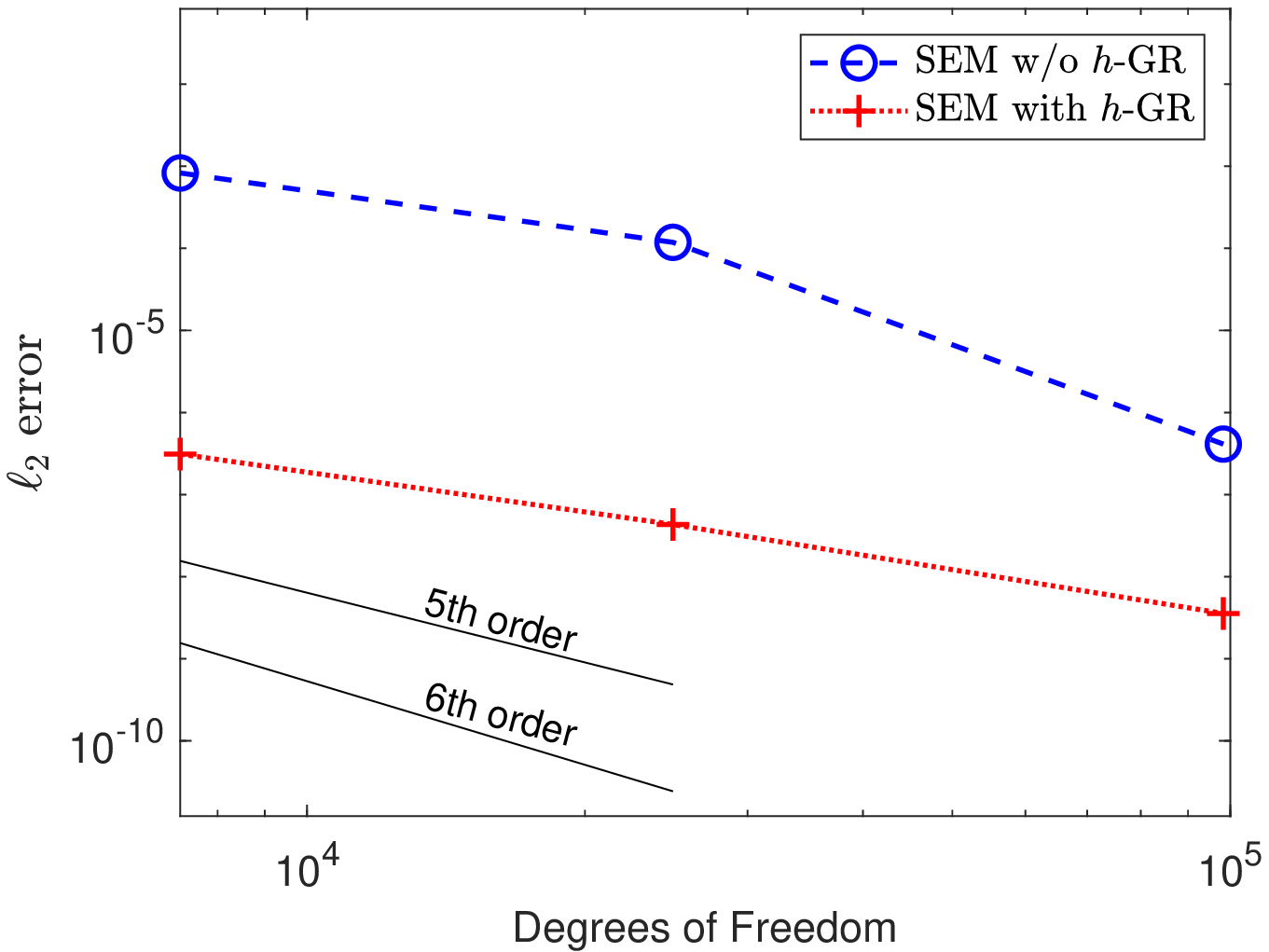}\\
{\scriptsize{}(c) Quartic}
\par\end{center}%
\end{minipage}
\par\end{centering}
\caption{\label{fig:CBR_flo_interior}Comparison of nodal errors of spectral
elements in the interior of the flower-hole domain with quadratic,
cubic, and quartic SEM for the convection-diffusion equation with
and without $h$-GR.}
\end{figure*}

To gain more insights, we plotted $\ell_{2}$-norm errors of the nodal
solutions of only the tensor-product spectral elements in the interior
of the domain for the flower-hole mesh in Figure~\ref{fig:CBR_flo_interior}.
It can be seen that curvature-based $h$-GR significantly improved
the accuracy of SEM. These results also confirmed our analysis in
section~\ref{subsec:Loss-and-Preservation} that the geometric errors
lead to significant pollution of the tensor-product elements in the
interior of the domain. These improvements contribute to the majority
of the improvements in Figure~\ref{fig:CBR_ell_iso}. Note that the
improvements from $h$-GR alone started to plateau as the meshes were
refined since the errors of near-boundary elements started to dominate,
which were further reduced by our post-processing technique as shown
in Figure~\ref{fig:domain_figures}.

\subsection{Effect of $p$-GR on Curved Neumann Boundaries\label{subsec:hp-refinement}}

\begin{figure*}
\begin{minipage}[t]{0.45\textwidth}%
\begin{center}
\includegraphics[width=1\textwidth]{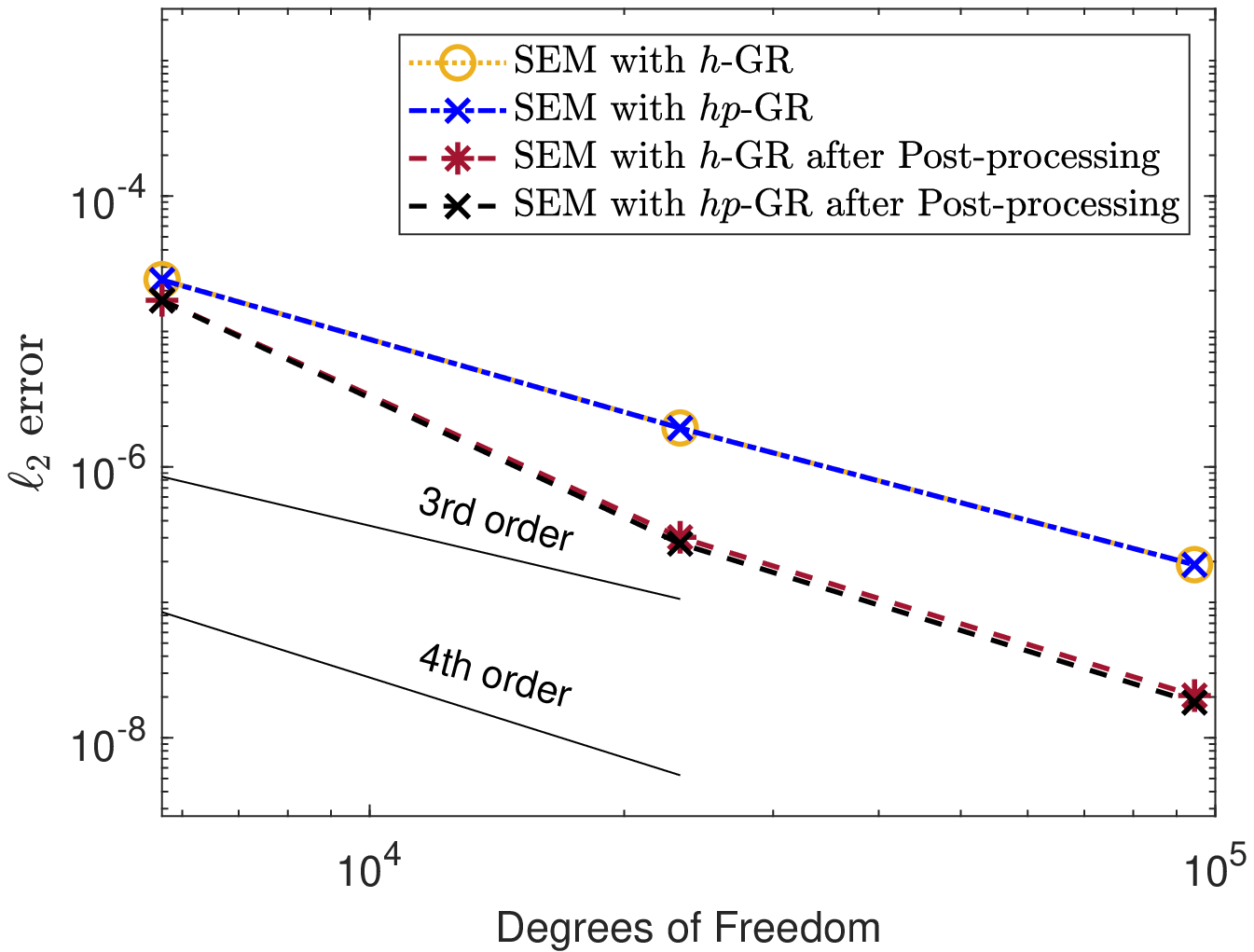}\\
{\scriptsize{}(a) Quadratic elliptical hole}
\par\end{center}%
\end{minipage} %
\begin{minipage}[t]{0.45\textwidth}%
\begin{center}
\includegraphics[width=1\textwidth]{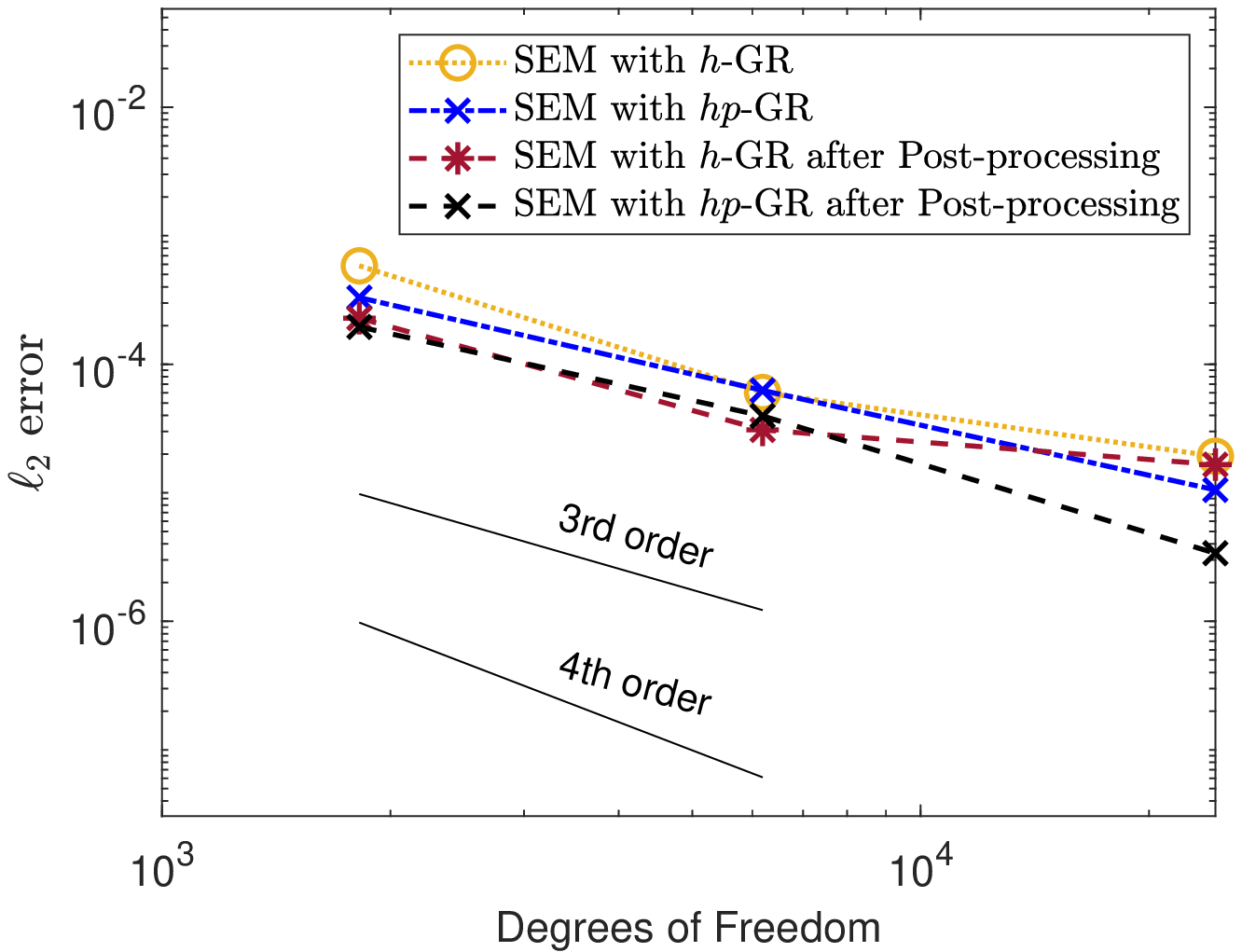}\\
{\scriptsize{}(b) Quadratic flower hole}
\par\end{center}%
\end{minipage}

\begin{minipage}[t]{0.45\textwidth}%
\begin{center}
\includegraphics[width=1\textwidth]{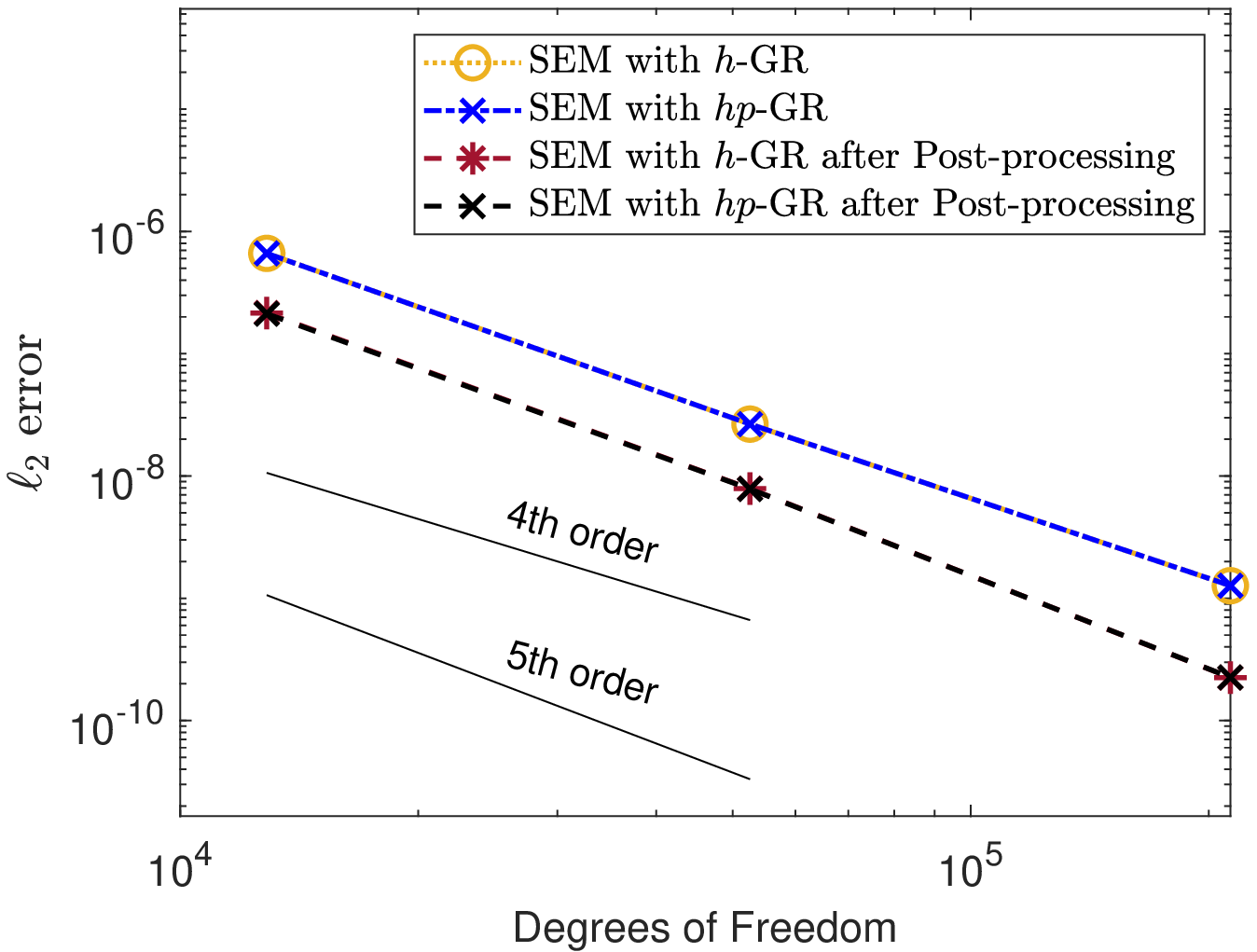}\\
{\scriptsize{}(c) Cubic elliptical hole}
\par\end{center}%
\end{minipage} %
\begin{minipage}[t]{0.45\textwidth}%
\begin{center}
\includegraphics[width=1\textwidth]{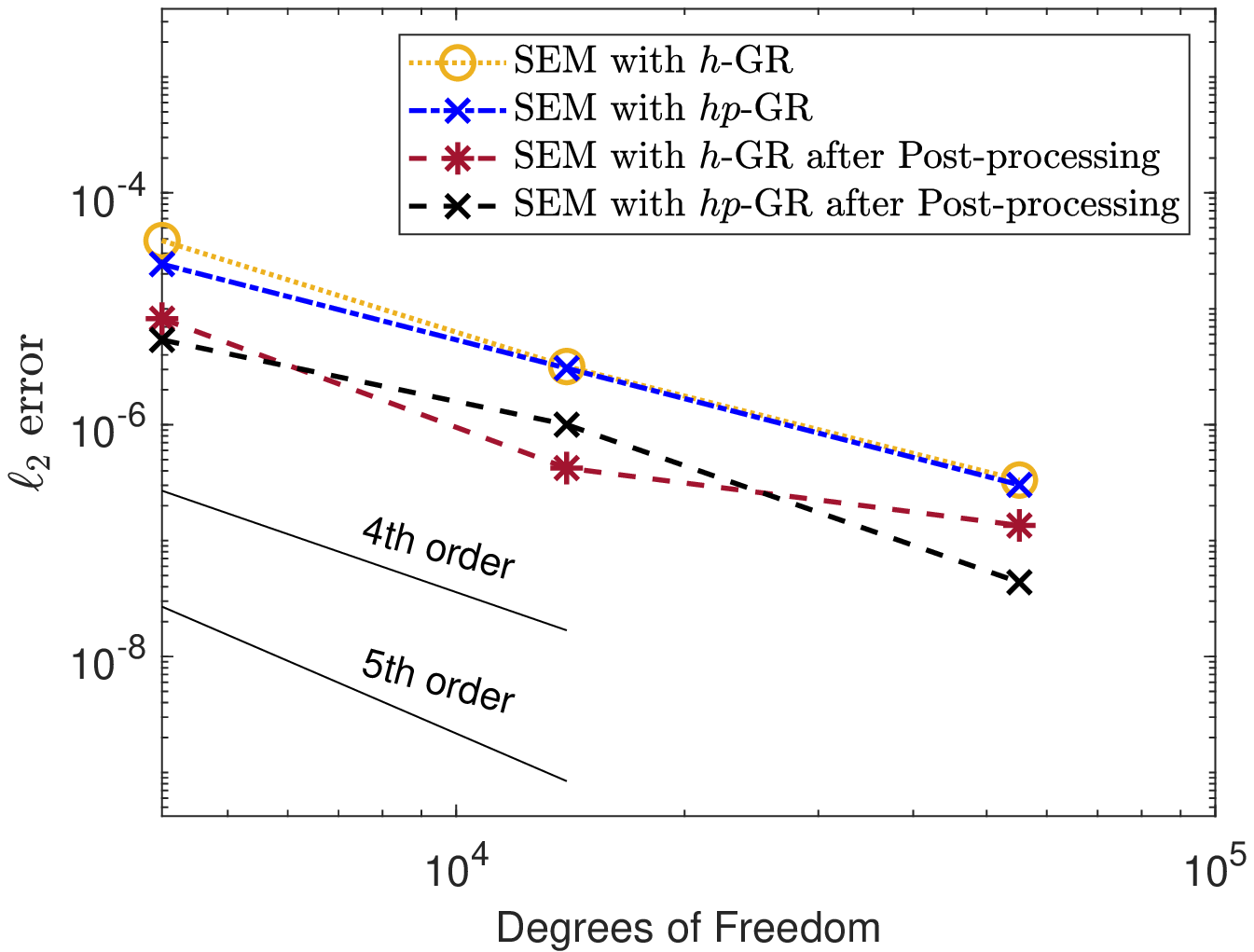}\\
{\scriptsize{}(d) Cubic flower hole}
\par\end{center}%
\end{minipage}

\begin{minipage}[t]{0.45\textwidth}%
\begin{center}
\includegraphics[width=1\textwidth]{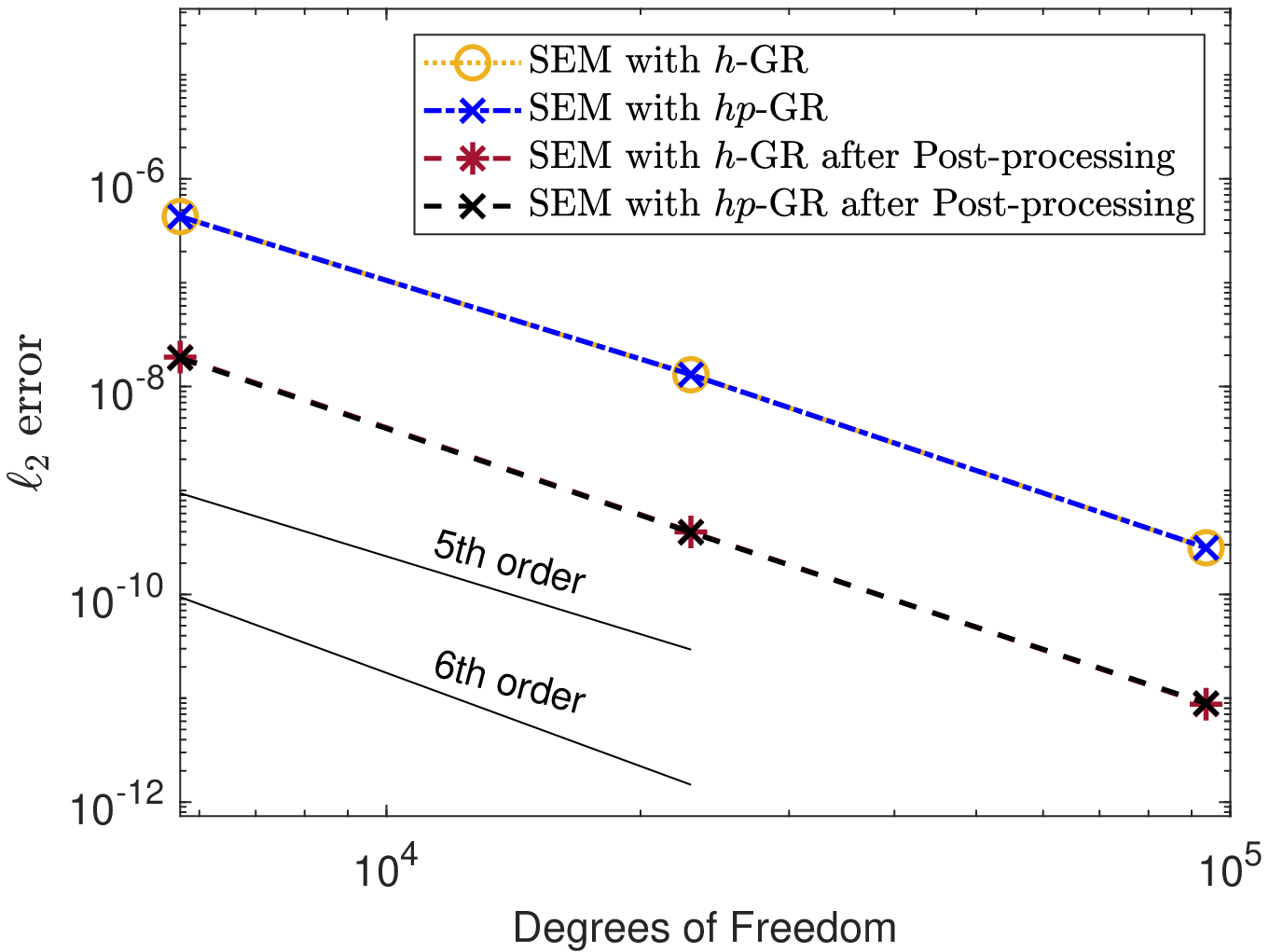}\\
{\scriptsize{}(e) Quartic elliptical hole}
\par\end{center}%
\end{minipage} %
\begin{minipage}[t]{0.45\textwidth}%
\begin{center}
\includegraphics[width=1\textwidth]{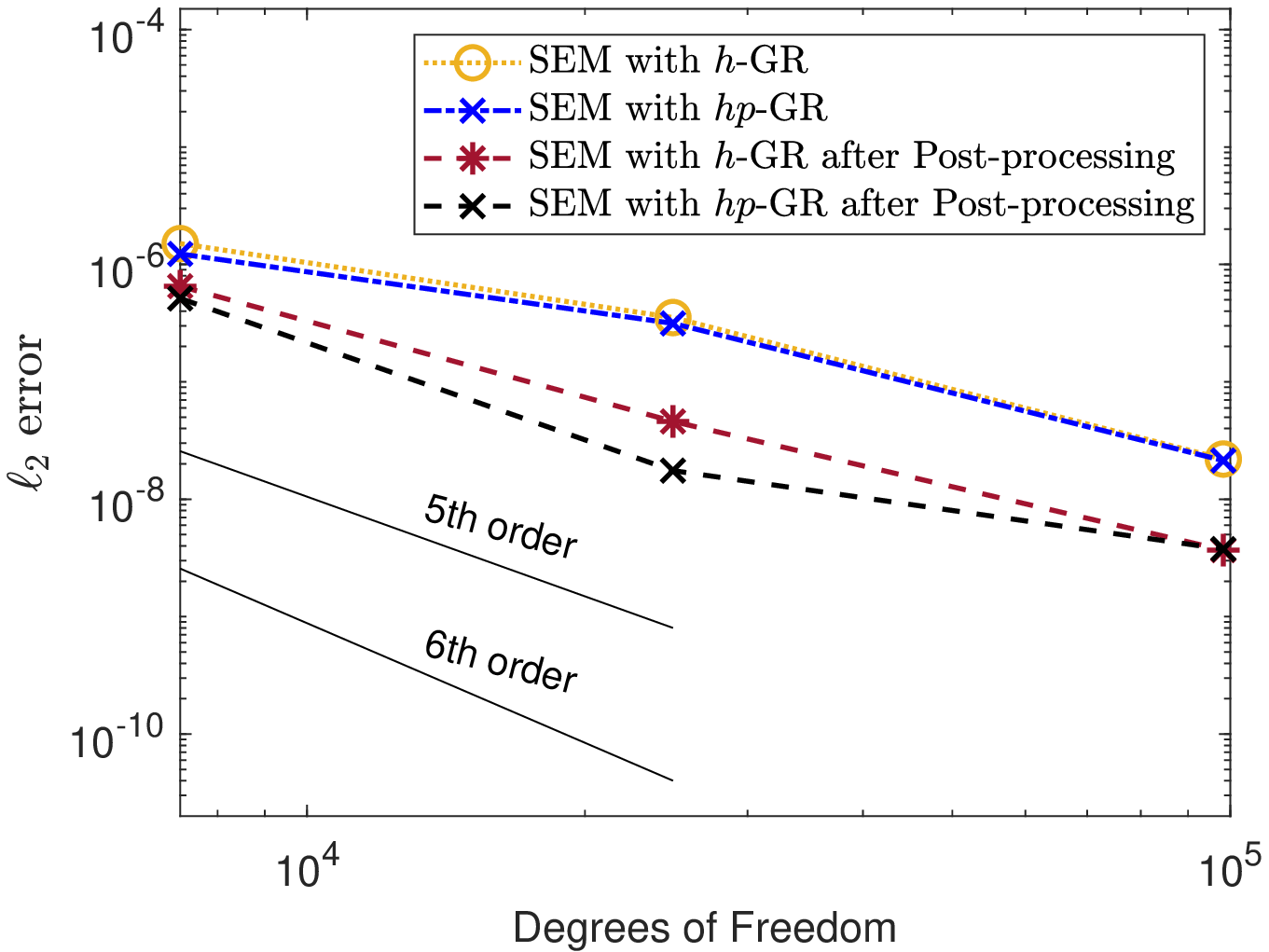}\\
{\scriptsize{}(f) Quartic flower hole}
\par\end{center}%
\end{minipage}

\caption{\label{fig:CBR_ell_isosuper}A comparison isoparametric elements versus
superparametric elements ($h$-GR vs $hp$-GR) for the convection-diffusion
equation on both the elliptical hole domain on the left side and the
flower hole domain on the right with Dirichlet boundary conditions
on the straight boundaries and Neumann boundary conditions on the
curved boundaries solved with quadratic, cubic, and quartic SEM as
well as the error after post-processing.}
\end{figure*}

Next, we study the effect of superparametric elements ($p$-GR) on
boundary elements. Similar to section~\ref{subsec:h-refinement},
we solved the convection-diffusion equation (\ref{eq:conv-diff})
on both the elliptical-hole and flower-hole meshes with Neumann and
Dirichlet boundary conditions on $\Gamma_{1}$ and $\Gamma_{2}$,
respectively. However, this time we compare the error using isoparametric
elements versus superparametric elements near the curved boundary
$\Gamma_{1}$ for meshes with $h$-GR. Figure~\ref{fig:CBR_ell_isosuper}
shows the $\ell_{2}$-norm errors for both $h$-GR and $hp$-GR using
SEM before and after post-processing. Here we omit the results with
$p$-GR only, since it may increase the error if the curved boundary
is under-resolved.

From the results, we observe that in a similar fashion to section~\ref{subsec:h-refinement},
$hp$-GR improved the geometric accuracy of normals, reducing the
error when using Neumann boundary conditions. However, the effect
of going from $h$-GR to $hp$-GR is much less dramatic using SEM
than from standard geometry to $h$-GR for the flower domain. The
more significant results become clear after post-processing, in which
superparametric elements can sometimes help recover superconvergence
for domains with high curvature as presented in section~\ref{subsec:Resolution-of-Curved}.
It has been theorized that isoparametric elements have sufficient
geometric accuracy and that superparametric elements do not improve
accuracy, in this regard. This result shows that when trying to recover
$\ell_{2}$ superconvergence that was lost due to geometric inaccuracies
near the boundary through post-processing, superparametric elements
may be very useful by reducing the interpolation errors.

\subsection{Accuracy Improvement for Curved Dirichlet Boundaries\label{subsec:DBC}}

In this subsection, we present the effect $h$-GR has using only Dirichlet
boundary conditions on all boundaries of the meshes. Each graph contains
the errors of SEM with and without $h$-GR as the errors after post-processing.
Similar to the results in sections~\ref{subsec:h-refinement} and
\ref{subsec:hp-refinement}, we observe that increasing the geometric
accuracy near a boundary with high curvature helps preserve the accuracy
of superconvergent SEM in the interior after post-processing. Dirichlet
boundary conditions, however, have some caveats in terms of how the
geometric accuracy impacts the $\ell_{2}$ error. In contrast with
nodes on Neumann boundary conditions, the nodes on a Dirichlet boundary
have known function values, so the difference in adjacent normal derivatives
does not play a role in the analysis. However, we do see that sufficient
geometric accuracy is required to maintain high-order accuracy for
interior elements after post-processing.

\begin{figure*}
\begin{minipage}[t]{0.45\textwidth}%
\begin{center}
\includegraphics[width=1\textwidth]{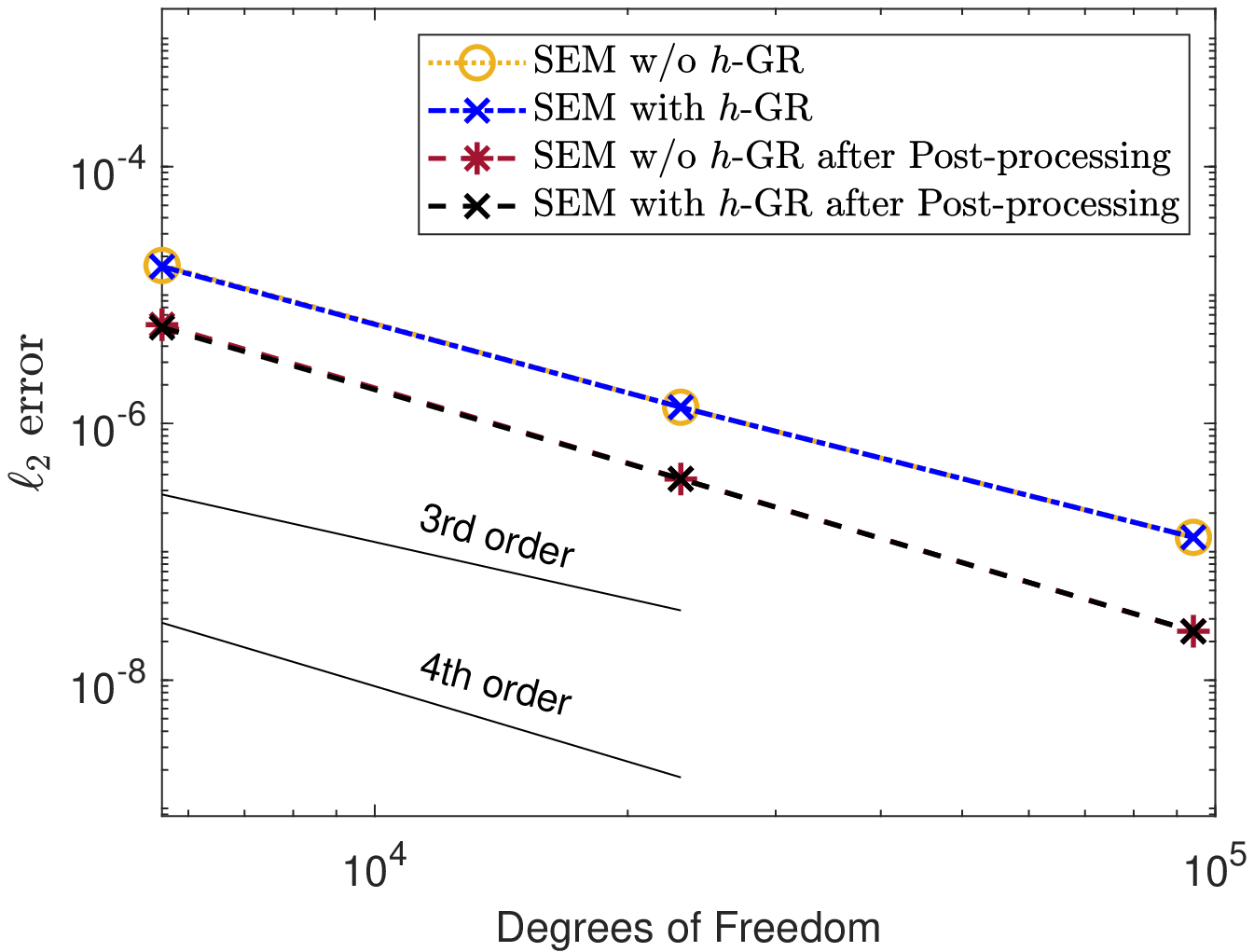}\\
{\scriptsize{}(a) Quadratic elliptical hole}
\par\end{center}%
\end{minipage} %
\begin{minipage}[t]{0.45\textwidth}%
\begin{center}
\includegraphics[width=1\textwidth]{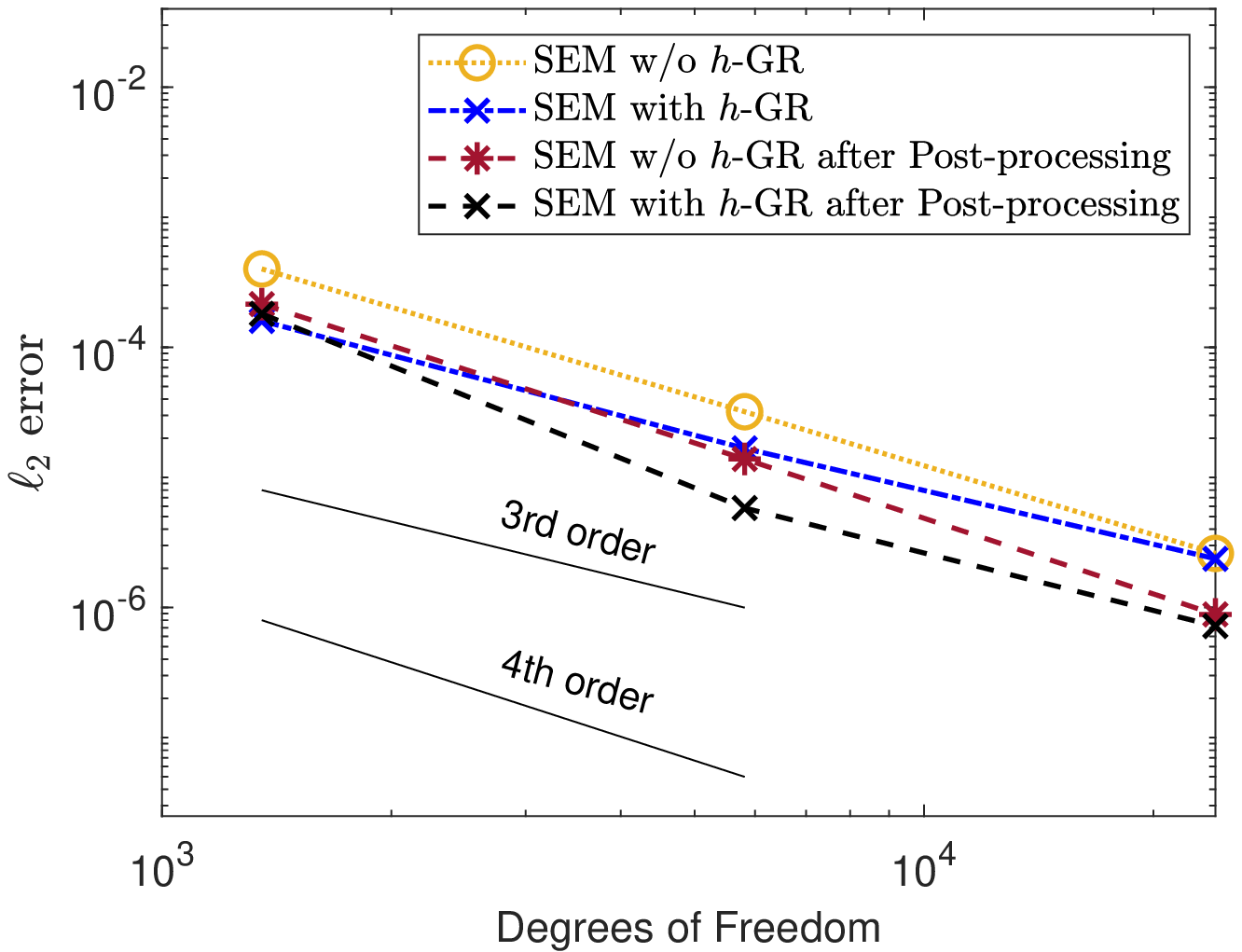}\\
{\scriptsize{}(b) Quadratic flower hole}
\par\end{center}%
\end{minipage}

\begin{minipage}[t]{0.45\textwidth}%
\begin{center}
\includegraphics[width=1\textwidth]{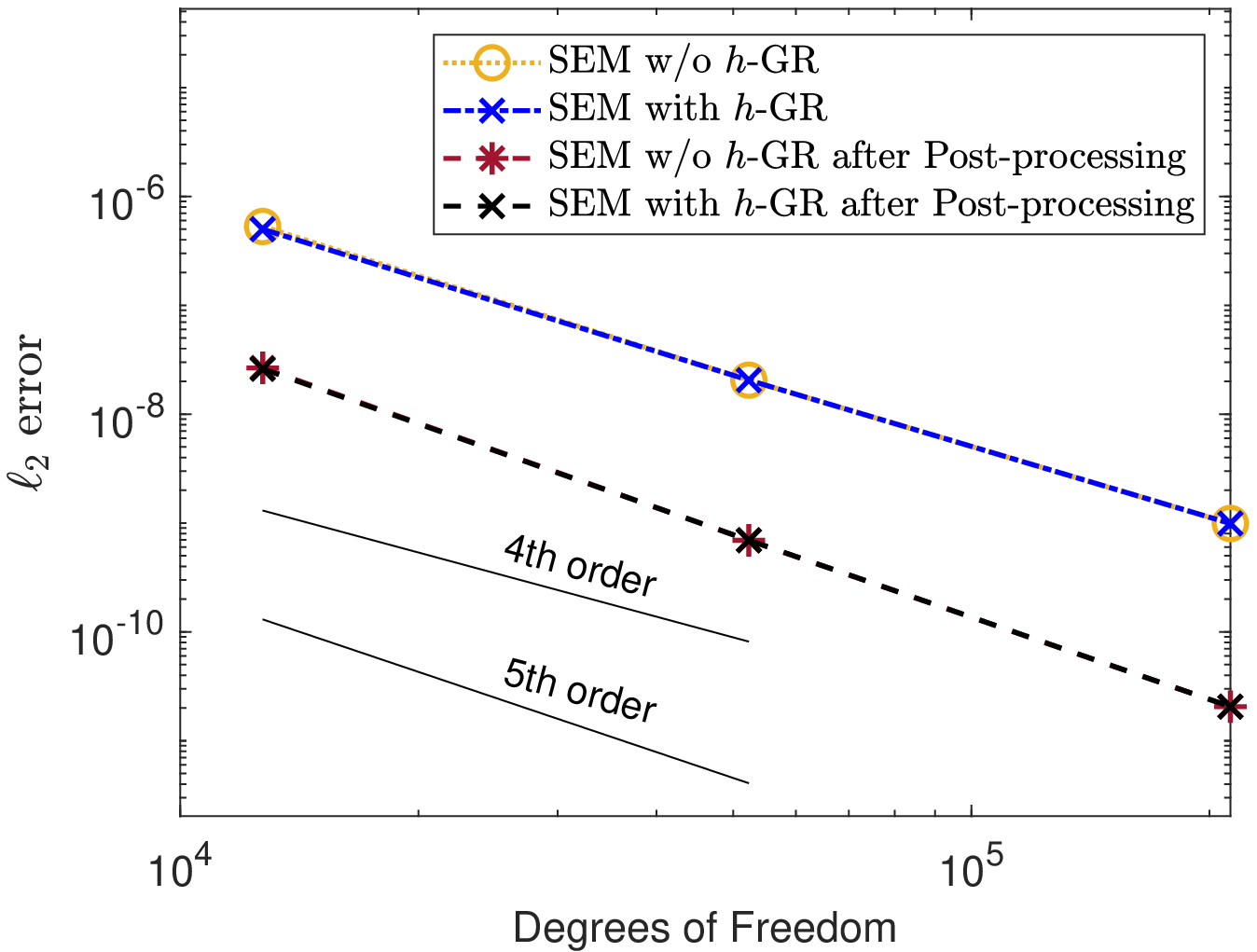}\\
{\scriptsize{}(c) Cubic elliptical hole}
\par\end{center}%
\end{minipage} %
\begin{minipage}[t]{0.45\textwidth}%
\begin{center}
\includegraphics[width=1\textwidth]{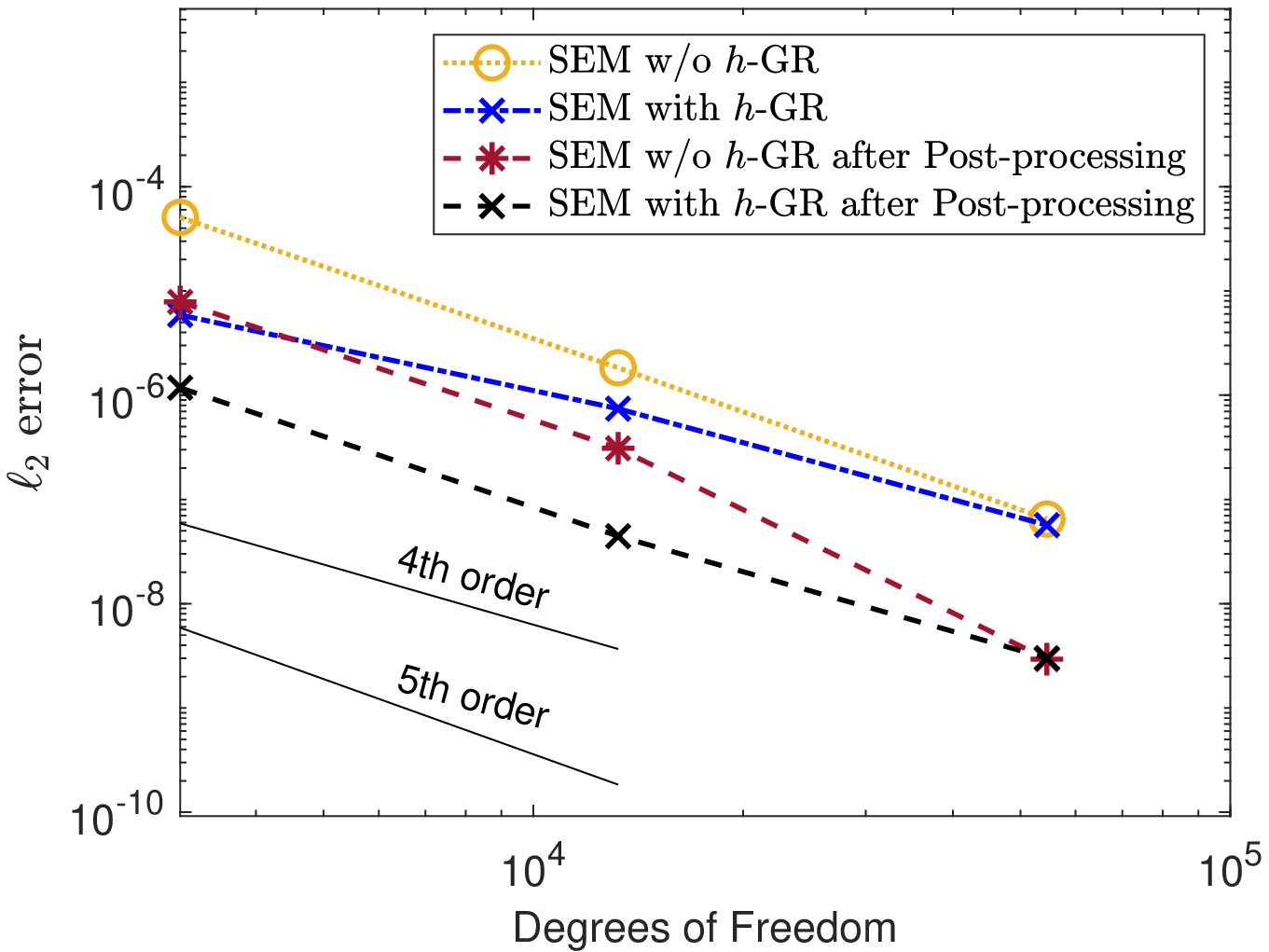}\\
{\scriptsize{}(d) Cubic flower hole}
\par\end{center}%
\end{minipage}

\begin{minipage}[t]{0.45\textwidth}%
\begin{center}
\includegraphics[width=1\textwidth]{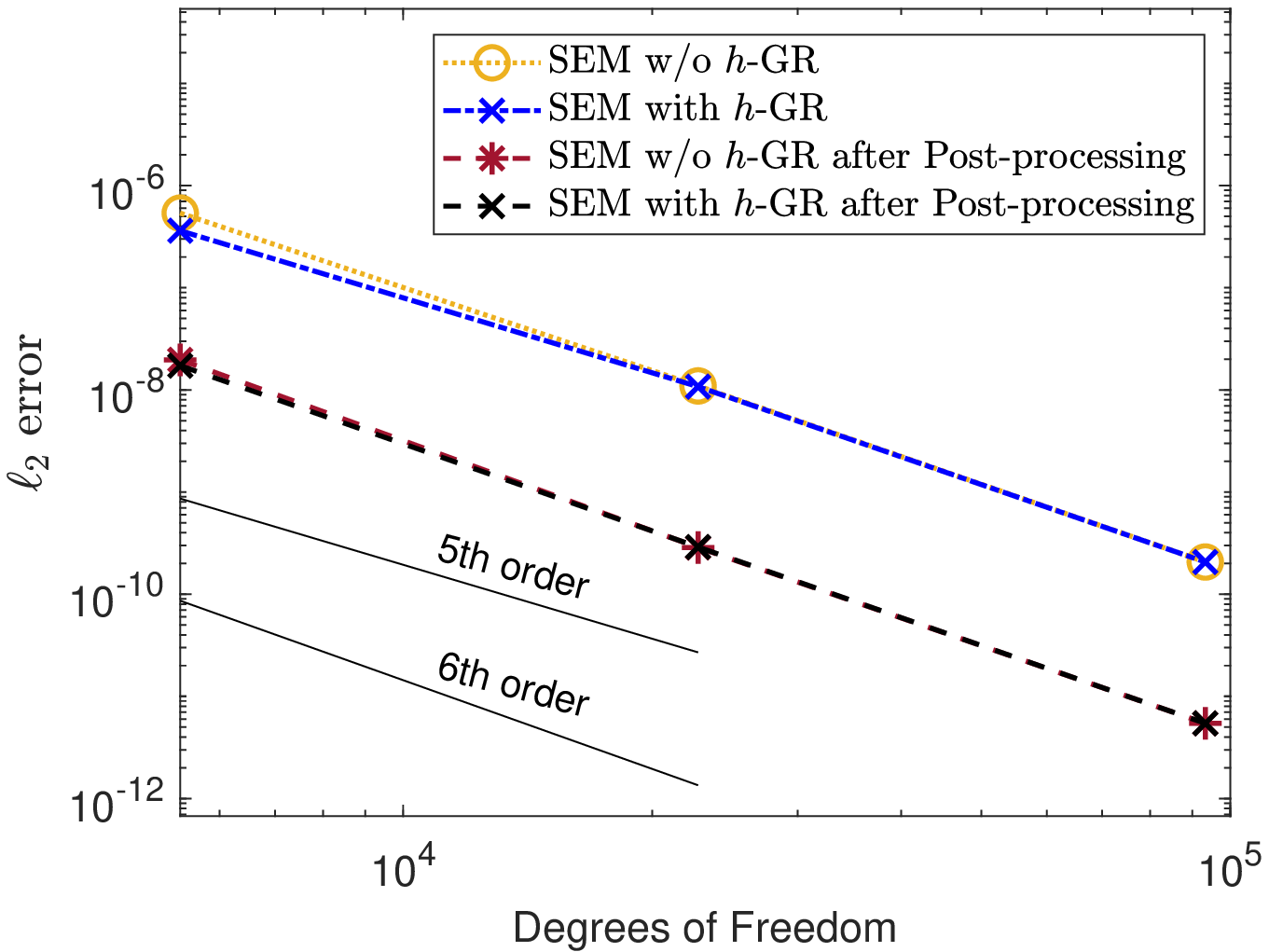}\\
{\scriptsize{}(e) Quartic elliptical hole}
\par\end{center}%
\end{minipage} %
\begin{minipage}[t]{0.45\textwidth}%
\begin{center}
\includegraphics[width=1\textwidth]{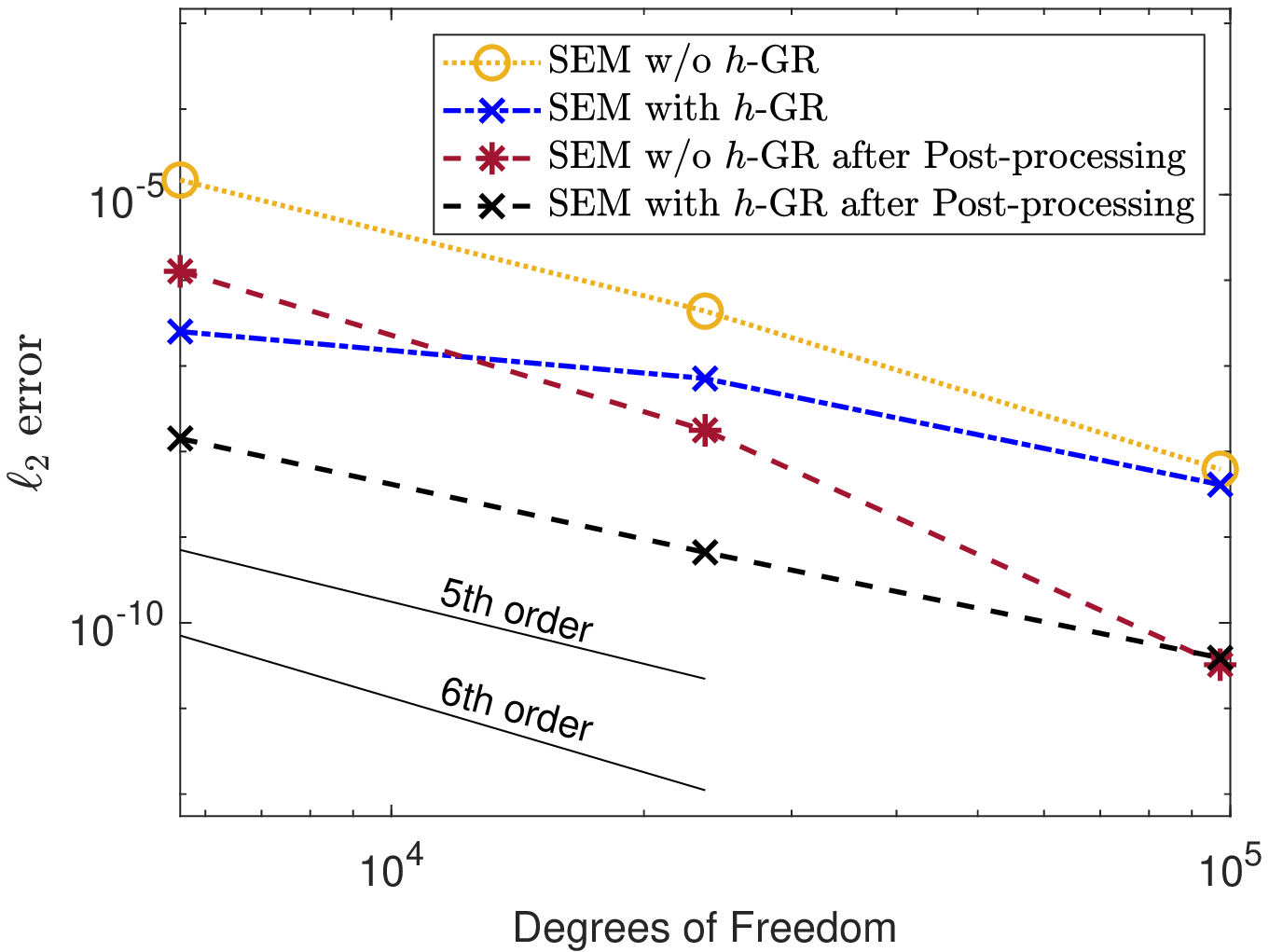}\\
{\scriptsize{}(f) Quartic flower hole}
\par\end{center}%
\end{minipage}

\caption{\label{fig:CBR_direll}A comparison of isoparametric elements versus
superparametric elements for the convection-diffusion equation on
both the elliptical hole domain on the left side and the flower hole
domain on the right with Dirichlet boundary conditions and quadratic,
cubic, and quartic SEM, with $h$-GR as well as the error after post-processing.}
\end{figure*}

As shown in the left side of Figure~\ref{fig:CBR_direll}, for the
meshes with the elliptical hole, $h$-GR does not offer any advantages
as standard meshes capture the true geometry close enough. For meshes
with regions of high curvature, $h$-GR offers more advantages since
meshes with $h$-GR better approximate the true geometry of the domain.
From the right side of Figure~\ref{fig:CBR_direll}, we can see that
for coarse meshes with $h$-GR the accuracy of both the SEM and the
post-processing solution are significantly better than the mesh without
$h$-GR. The standard mesh can better approximate exact geometry for
finer meshes, so the error is similar to using $h$-GR.

\subsection{Impact of Mesh Quality\label{subsec:meshquality}}

\begin{figure*}
\begin{centering}
\begin{minipage}[t]{0.5\textwidth}%
\begin{center}
\includegraphics[width=1\textwidth]{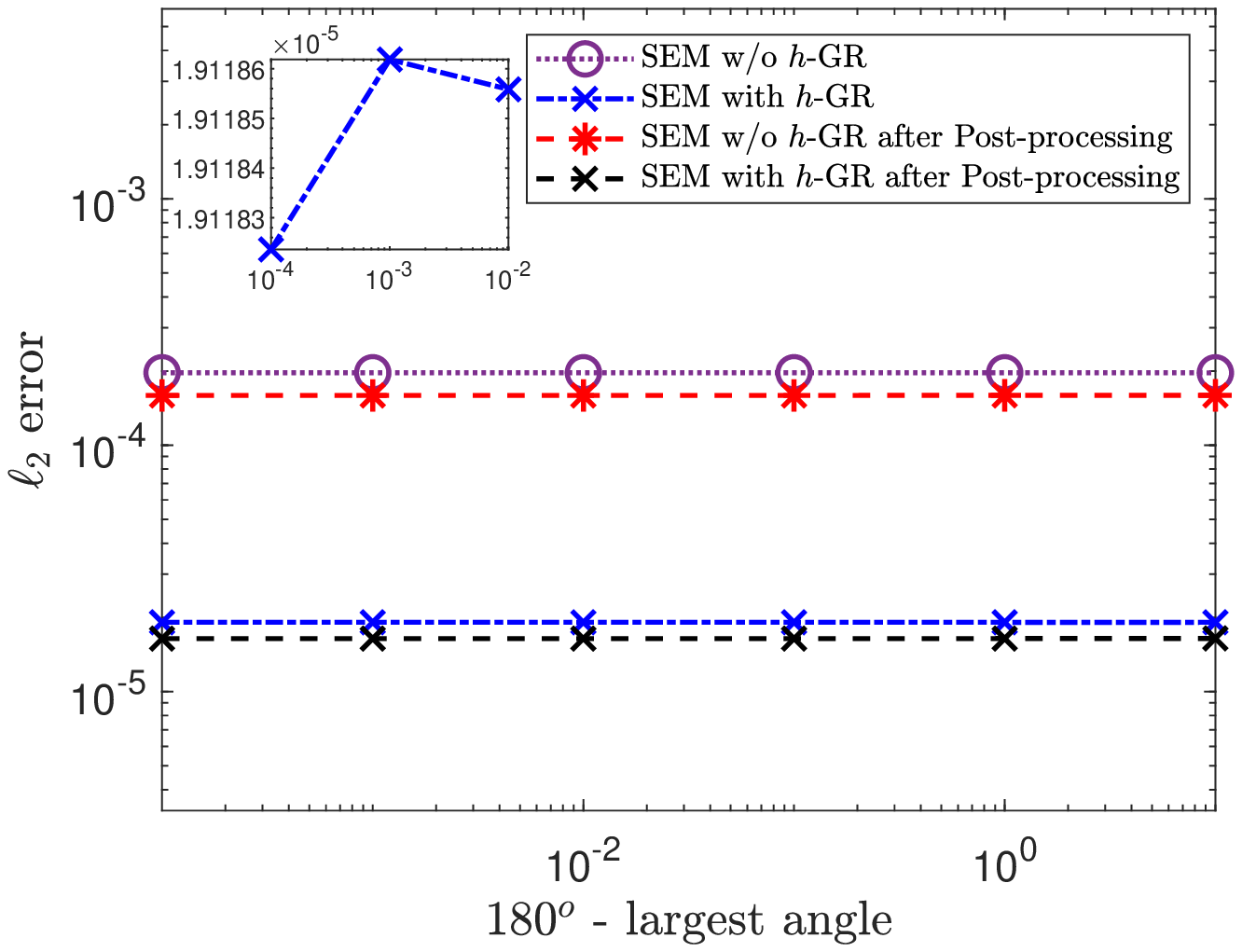}\\
{\scriptsize{}(a) Quadratic}
\par\end{center}%
\end{minipage}%
\begin{minipage}[t]{0.5\textwidth}%
\begin{center}
\includegraphics[width=1\textwidth]{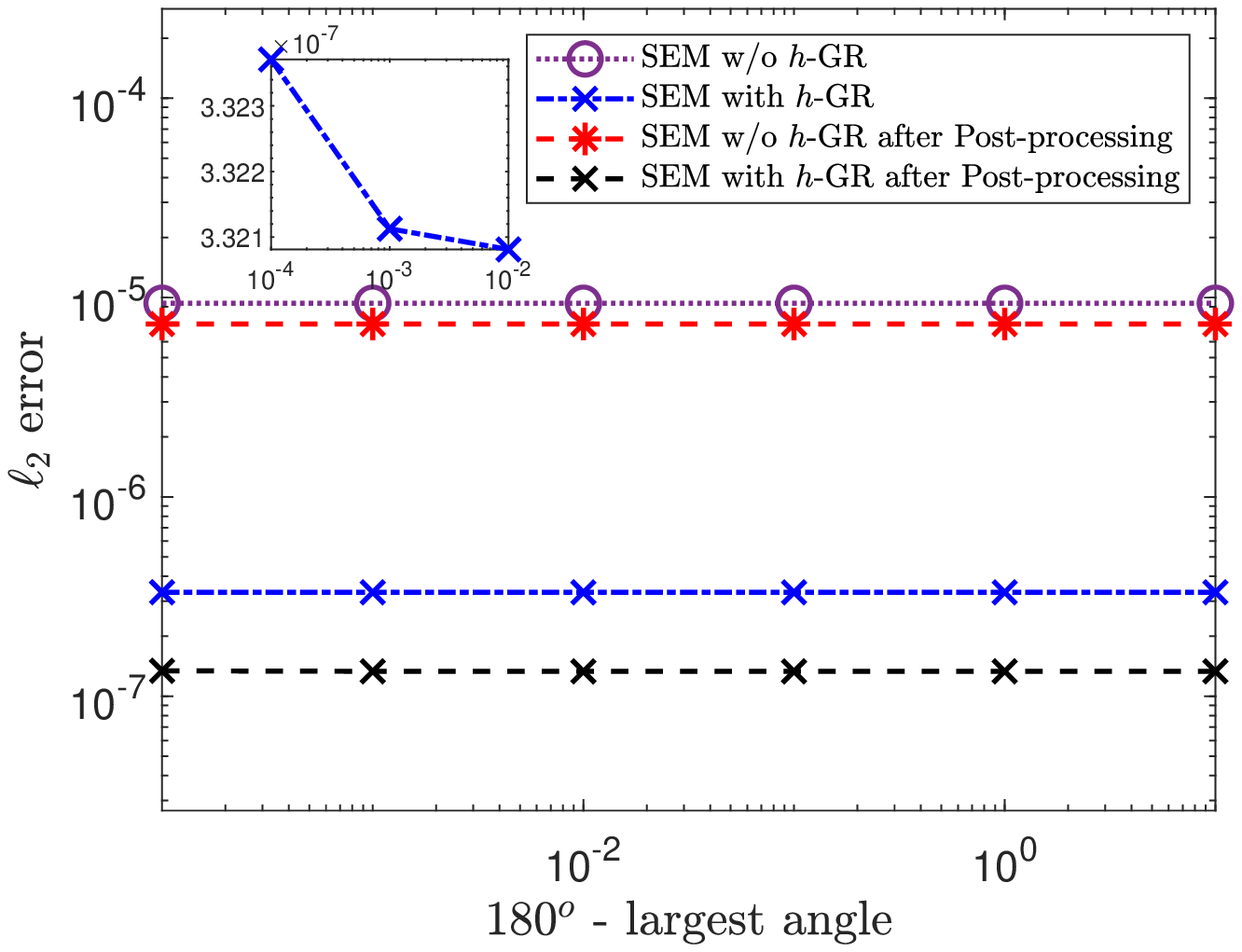}\\
{\scriptsize{}(b) Cubic}
\par\end{center}%
\end{minipage}\\
\begin{minipage}[t]{0.5\textwidth}%
\begin{center}
\includegraphics[width=1\textwidth]{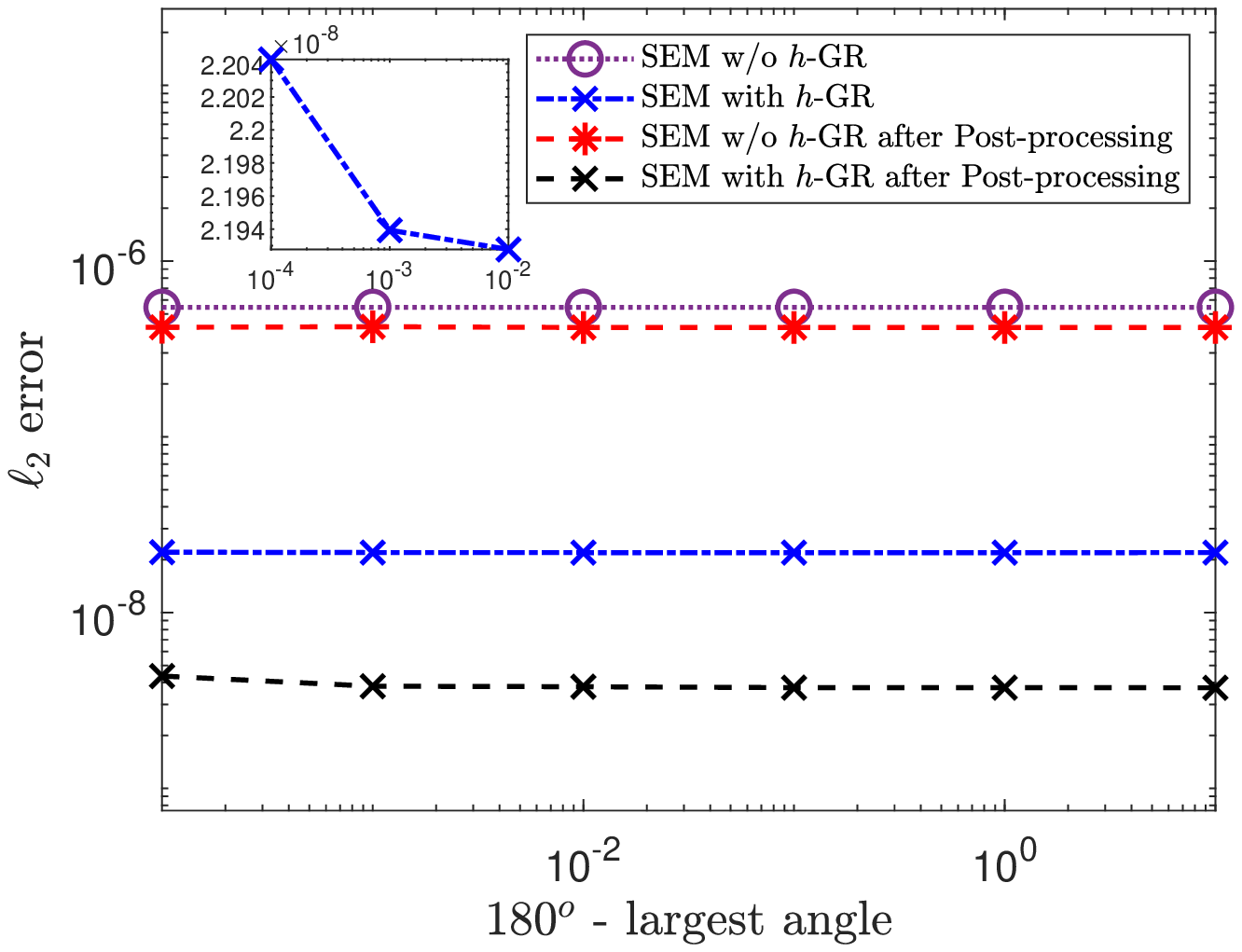}\\
{\scriptsize{}(c) Quartic}
\par\end{center}%
\end{minipage}
\par\end{centering}
\caption{\label{fig:CBRmeshqual}A mesh quality study which compares the errors
of SEM and post-processing with and without $h$-GR. The horizontal
axis shows the difference of $180\lyxmathsym{\protect\textdegree}$
minus the largest angle in the mesh. The inset is a zoomed in version
of the first three SEM errors.}
\end{figure*}

Finally, we present the effect that mesh quality has on the $l_{2}$-error
of meshes with and without $h$-GR, and the improvement of post-processing
on the same meshes. Since we wanted to look at the effect geometric
accuracy has on high curvatures, we will only use domain (b) from
Figure~\ref{fig:domain_figures}. We manipulated the mesh quality
by distorting a good quality mesh by moving several nodes such that
the maximum angle among all elements is some specified degree $\theta.$
Specifically, we consider a series of meshes where the maximum angle
is $180$ degrees minus $[10^{-4},10{}^{-3},10{}^{-2},10{}^{-1},10^{0},10{}^{1}]$
degrees. The results are shown in Figure~\ref{fig:CBRmeshqual},
which plots the $l_{2}$-error on the $y$-axis and $180\text{°}$
minus the maximum angle on the $x$-axis. The figure inset shows the
first three error values of SEM.

First, we note that post-processing can improve the accuracy more
as the order increases. We see this specifically as while post-processing
meshes with $h$-GR improves the accuracy a little for the quadratic
mesh, the quartic mesh enjoyed more than a ten-fold improvement in
nodal error. This result is due to the improvement of geometric accuracy
allowing for recovery of superconvergence as outlined in subsection~\ref{subsec:Loss-and-Preservation}.
We also want to point out that meshes without $h$-GR do not improve
much after post-processing, so superconvergence was lost. In Figure~\ref{fig:CBRmeshqual},
we also show that $h$-GR and post-processing can reduce error significantly
more than mesh quality optimization, provided that the maximum angle
is less than $179.999\text{°}$, a relatively relaxed mesh generation
requirement. These results suggest that $h$-GR and post-processing
may greatly improve angle quality requirements for mesh generation
and could reduce the meshing bottleneck significantly. Using AES-FEM
near the boundary relaxes the dependency on element shapes. It will
also ease the generalization to 3D, for which tetrahedral meshes are
more prone to having poor-shaped elements (such as slivers) near boundaries.

\section{\label{sec:Conclusions}Conclusions}

In this paper, we introduced a novel approach to improve the accuracy
of SEM over curved domains. Our approach utilizes mixed meshes with
superconvergent tensor-product spectral elements in the interior of
the domain and non-tensor-product elements (such as simplicial elements)
near the boundary. We showed that $h$- and $p$-GR of near-boundary
elements could substantially reduce the pollution errors to the spectral
elements in the interior and preserve the superconvergence property
of spectral elements. In addition, by utilizing AES-FEM to post-process
the solutions near the boundary, we can make the accuracy of the non-tensor-product
elements match those of the spectral elements. To the best of our
knowledge, this work is the first to demonstrate superconvergence
of SEM over curved domains. We showed that our proposed techniques
could improve the accuracy of the overall solution by one to two orders
of magnitude. At the same time, the element shapes have a substantially
smaller impact on the accuracy. As a result, our proposed approach
can substantially relax the mesh-quality requirement and, in turn,
alleviate the burden of mesh generation for higher-order methods,
especially near boundaries. We presented numerical results in 2D as
proof of concept. Future work includes the implementation in 3D and
optimization of the implementation.

\section*{Declarations}

\paragraph{Conflict of interest }

The authors have no relevant financial or non-financial interests
to disclose.

\bibliographystyle{abbrv}
\bibliography{ref/refs}

\end{document}